\documentclass[aps,pre,twocolumn,superscriptaddress,showpacs,showkeys ]{revtex4-1}

\usepackage{graphicx}                   % for insert file eps
\usepackage{hyperref}                   % for hyper refrences
\usepackage{amsmath,amssymb}            % for insert extra math
\usepackage{epstopdf}
\usepackage{placeins}

\usepackage{color}
\definecolor{orange}{rgb}{1,0.5,0}
\definecolor{brown}{rgb}{0.65, 0.16, 0.16}
\definecolor{phlox}{rgb}{0.87, 0.0, 1.0}

\graphicspath{{figs/}}					% the path of image folder 

\begin{document}
	
	\title{Homology Groups of Embedded Fractional Brownian Motion}
	%\title{Topological Distribution of Embedded Fractional Brownian Motion}

	\author{H. Masoomy}
	\email{hoseingmasoomy@gmail.com}
	\affiliation{Department of Physics, Shahid Beheshti University,  1983969411, Tehran, Iran}
	
	\author{S. Tajik}
	\email{samintjk@gmail.com}
	\affiliation{Department of Physics, Brock University, St. Catharines, ON L2S 3A1, Canada}
	
	\author{S. M. S. Movahed}
	\email{s.movahed@sbu.ac.ir}
	\affiliation{Department of Physics, Shahid Beheshti University,  1983969411, Tehran, Iran}

	\begin{abstract}
	A well-known class of non-stationary self-similar time series is the fractional Brownian motion (fBm) considered to model ubiquitous stochastic processes in nature. In this paper, we study the homology groups of high-dimensional point cloud data (PCD) constructed from synthetic fBm data.  We covert the simulated fBm series to a PCD, a subset of unit $D$-dimensional cube, employing the time delay embedding method for a given embedding dimension and a time-delay parameter. In the context of persistent homology (PH), we compute topological measures for embedded PCD as a function of associated Hurst exponent, $H$, for different embedding dimensions, time-delays and amount of irregularity existed in the dataset in various scales. Our results show that for a regular synthetic fBm, the higher value of the embedding dimension leads to increasing the $H$-dependency of topological measures based on zeroth and first homology groups. To achieve a reliable classification of fBm, we should consider the small value of time-delay irrespective of the irregularity presented in the data. More interestingly, the value of scale for which the PCD to be path-connected and the post-loopless regime scale are more robust concerning irregularity for distinguishing the fBm signal. Such robustness becomes less for the higher value of embedding dimension. 
	\end{abstract}

	
	\maketitle

	\section{Introduction}
	
Data generation from a wide range of experiments and observations leads to great opportunities to use modeling in the context of complex systems. Assigning a \textit{shape} to various types of data and accordingly computing their persistence in terms of so-called resolution has recently become the focus of many studies. Furthermore, the spatio-temporal distribution (shape) of datasets has been beneficial in various branches of science ~\cite{vuran2004spatio, diggle2006spatio, chaplain2001spatio}. Among different criteria such as linear algebra~\cite{strang1993introduction,defranza2015introduction} and common statistical properties ~\cite{jaynes1957information,goulden1939methods, mandel2012statistical,albert2002statistical,costa2007characterization,barthelemy2011spatial,wu2015emergent,boguna2021network}, the geometrical and topological features has provided complementary evaluations ~\cite{adler1981geometry,carlsson2009topology}.
Depending on what type of information is needed and based on limitations in computational resources, one- and/or $n$-point statistics of geometrical and topological properties \cite{rice44a,rice44b,rice1954selected,Bardeen:1985tr,matsubara03,Pogosyan:2008jb,Gay:2011wz,matsubara2020large,vafaei2021clustering}, and complex network based analysis \cite[and references therein]{albert2002statistical,costa2007characterization,barthelemy2011spatial,wu2015emergent,boguna2021network}, can be taken into account.  

Inspired by geometrical fractals~\cite{mandelbrot1982fractal}, the notion of self-similarity and self-affinity is employed to quantify the statistical properties of different fields of generally $(N+D)$-dimension\cite{adler1981geometry} \footnote{A typical $(N+D)$-dimensional stochastic field (process), is a measurable mapping from probability space into a $\mathbb{R}^{N}$-valued function on $\mathbb{R}^{D}$ Euclidean space, in which the labels $N$ and $D$ are respectively devoted to $N$ dependent and $D$ independent parameters, generating a $(N+D)$-dimensional random field or a stochastic process}. A pioneer method to quantify a self-similar process based on Hurst exponent \cite{hurst1951long} which is known as fractional Brownian motion (fBm) (its associated increment is called fractional Gaussian noise (fGn)) \cite{mandelbrot1968fractional}, is called Multifractal Detrended Fluctuation Analysis (MFDFA) \cite{peng1994mosaic,Peng95,Kantelhardt} (see also \cite[and references therein]{jiang2019multifractal,eghdami2018multifractal}). Besides the noted methods for classification, prediction, and examination of datasets, one can focus on topological features of data to study it under the banner of Topological Data Analysis (TDA)~\cite{dey2022computational, edelsbrunner2022computational, zomorodian2005topology, ghrist2008barcodes, wasserman2018topological} and accordingly construct a powerful algebra-topological-based approach~\cite{nakahara2003geometry, munkres2018elements, hatcher2005algebraic}. To this end, Persistent Homology (PH) as a shape-based tool that examines the evolution of global features of the data related to topological invariants is utilized~\cite{otter2017roadmap, speidel2018topological, maletic2016persistent, donato2016persistent, pereira2015persistent, masoomy2021topological}. 
	
There are many methods to construct higher dimensional sets from a $(1+1)$-dimensional series and then implement TDA. The visibility graph  method~\cite{lacasa2008time}, correlation network method~\cite{yang2008complex} and time delay embedding (TDE) method~\cite{takens1981detecting, packard1980geometry} have been used as the pre-processing part on the input data. The visibility graph and correlation network methods constructe a complex network from the time series by considering the visibility condition between data points of underlying data and correlation between sub-time series of the dataset, respectively, while the TDE method maps the time series into a $D$-dimensional point cloud data (PCD) by a time-delay parameter, $\tau$. There are also some undirected methods for converting time series to the complex network, in which the reconstructed PCD (by TDE method) is converted to the network. The recurrence network method~\cite{gao2009flow, marwan2009complex, donner2010recurrence} and $k$-nearest neighbor (kNN) network method~\cite{shimada2008analysis} are the examples. Briefly, in the recurrence network method, the connections are considered based on the proximity of the embedded data points (state vectors), while in the kNN method all state vectors (nodes) are connected with their $k$ nearest neighbors.\\

	Quantifying the scaling exponent of fBm and fGn has been done by various methods \cite{peng1994mosaic,Peng95,Kantelhardt,jiang2019multifractal,eghdami2018multifractal}, but the specious effect of complicated trends, the impact of non-stationarities, the finite size effect and irregularities have remained as the main challenges in the mentioned algorithms  \cite{hu2001effect,trend2,nagarajan2005minimizing,ma2010effect,eghdami2018multifractal}. Subsequently, some efforts have been devoted to removing or at least tuning down the above discrepancies from different points of view ranging from network analysis \cite{lacasa2009visibility, ni2009degree,xie2011horizontal,masoomy2021visibility},  reconstructed phase space of fBm series by TDE method in terms of recurrence network analysis \cite{liu2014topological} by taking into account the false nearest neighbors algorithm \cite{kennel1992determining} and mutual information  method~\cite{fraser1986independent}, recurrence network of fBms \cite{zou2015analyzing}, to performing TDA on the weighted natural visibility graph constructed from fGn series \cite{masoomy2021persistent}. 
	
	The main purpose of this work is to figure out how the reconstructed phase space (embedded PCD) of the fBm time series changes by varying the intrinsic parameter of the signal, i.e. the Hurst exponent, $H$, and the algorithmic parameters of the TDE method, i.e. the embedding dimension, $D$, and the time-delay, $\tau$. Furthermore, the influence of irregularity adjusted by a parameter, $q$, the fraction of the number of missing data points, $T_{\rm missing}$, in the time series to its total number of data points, $T$, as a measure for the amount of irregularity of underlying dataset will be examined in this work. More precisely, we utilize the PH method to get deep insight into how the state vectors constructed by the TDE method are distributed from a topological viewpoint. One of the advantages of this idea is that we can study the global properties of the phase space measured by the population of the homology groups of the weighted topological space (called weighted simplicial complex) mapped from the corresponding reconstructed phase space. And more importantly, we can also capture the evolution of these topological invariants by varying the proximity parameter (threshold), $\epsilon$, continuously. Subsequently, we examined the behavior of the topological measures computed from the $d$th persistent diagram, $\mathcal{D}_{d}$, and $d$th Betti number, $\beta_{d}$, as a function of the proximity parameter (called $d$th Betti curve) for $d=0,1$. We also aim to find the optimal choice for the parameters $D$ and $\tau$ to distinguish the fBms of various $H$ and also the results would be robust against the irregularity. %(finding the best choice for $D$ and $\tau$ such that the computed measures classify the time series by their Hurst exponent and also have sufficient robustness to irregularity).
		
	Our results indicate that the computed topological measures of embedded fBm are sensitive to the value of the Hurst exponent. Generally, the statistics of $d$th homology classes have strong Hurst-dependency. The evolution of  $d$-dimensional topological holes ($d$-holes) occurs in low (high) scales for fBms with $H\gtrsim0.5$ ($H\lesssim0.5$). This $H$-dependency grows by increasing the dimension of constructed PCD, i.e. for a good estimation for $H$ one can deal with high-dimensional PCD. The situation is more noticeable when the time series contains some irregularity. For a signal with a high value of irregularity,  the proposed measures become more $D$-dependent and the accuracy of estimating the Hurst exponent decays for high-dimensional PCD. In other words, the topological features extracted from low-dimensional PCDs have the minimum $q$-dependency, i.e. the features of 2-dimensional reconstructed phase space are more robust than the $D$-dimensional PCD for $D>2$ in the presence of irregularity, suggesting the best value $D=2$ for the estimation of Hurst exponent of irregular fBms. 
	
	The rest of this paper is organized as follows: In the next section, our methodology and pipeline for analysis of the fBm signal are introduced. The numerical results of synthetic fBm time series, via reconstructed phase-space distribution of the state vectors from the topological viewpoint, are given in section \ref{results}. Summary and concluding remarks will be presented in section \ref{summary}.

	\begin{figure}
		\begin{center}
			\includegraphics[width=0.45\textwidth]{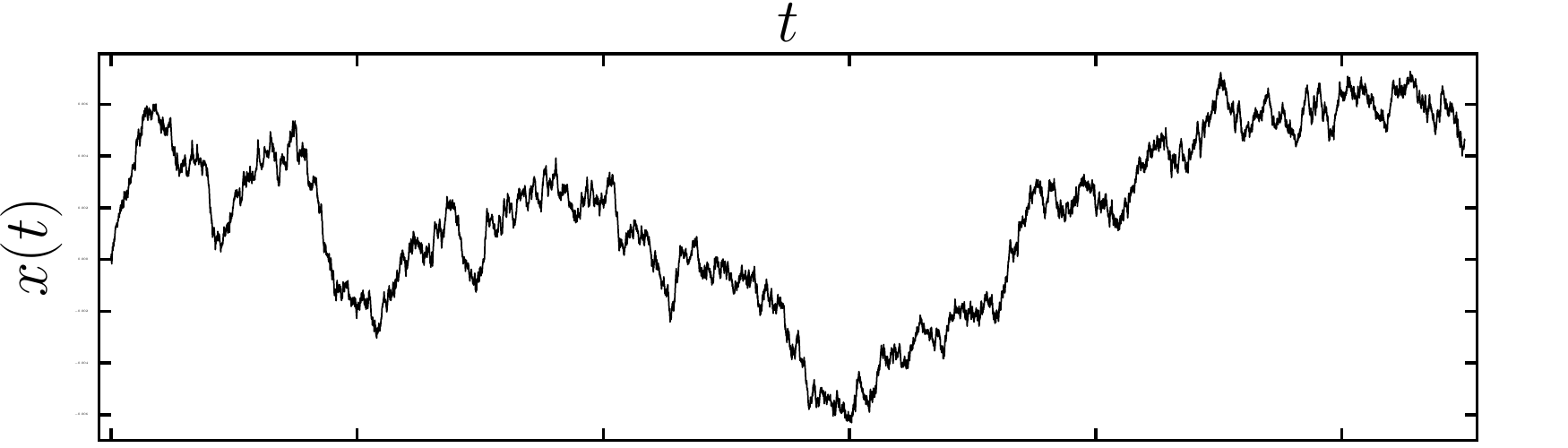}
			\includegraphics[width=0.45\textwidth]{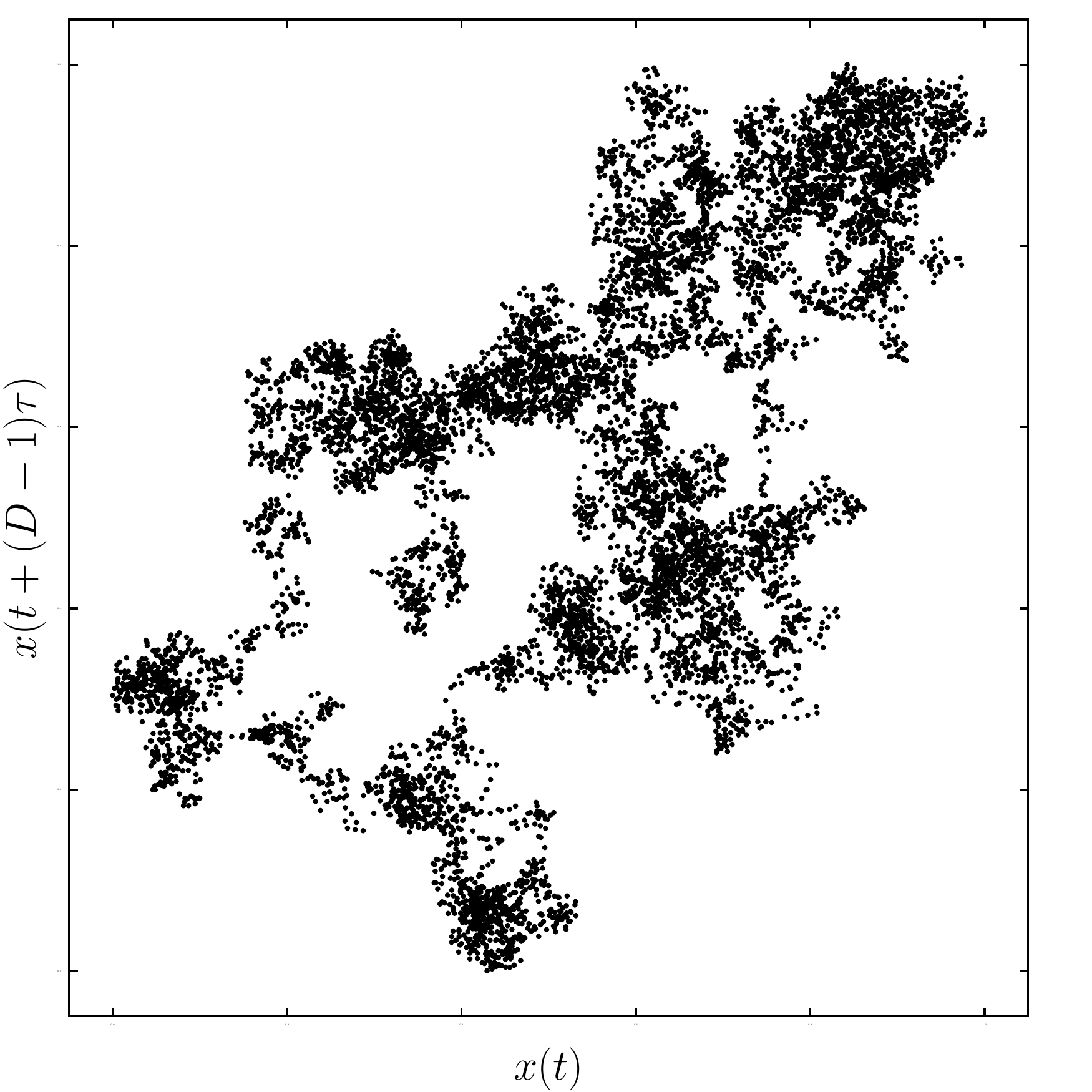}
			\caption{A synthetic regular ($q=0$) fBm time series of length $T=10^{4}$ with $H=0.5$ (top panel) and associated phase space of size $N=9\times10^{3}$ reconstructed by the TDE method for embedding dimension $D=2$ and time-delay $\tau=10^{3}$ (bottom panel).}
			\label{fig:tde}
			
		\end{center}
	\end{figure}
	\begin{figure}
		\includegraphics[width=0.45\textwidth]{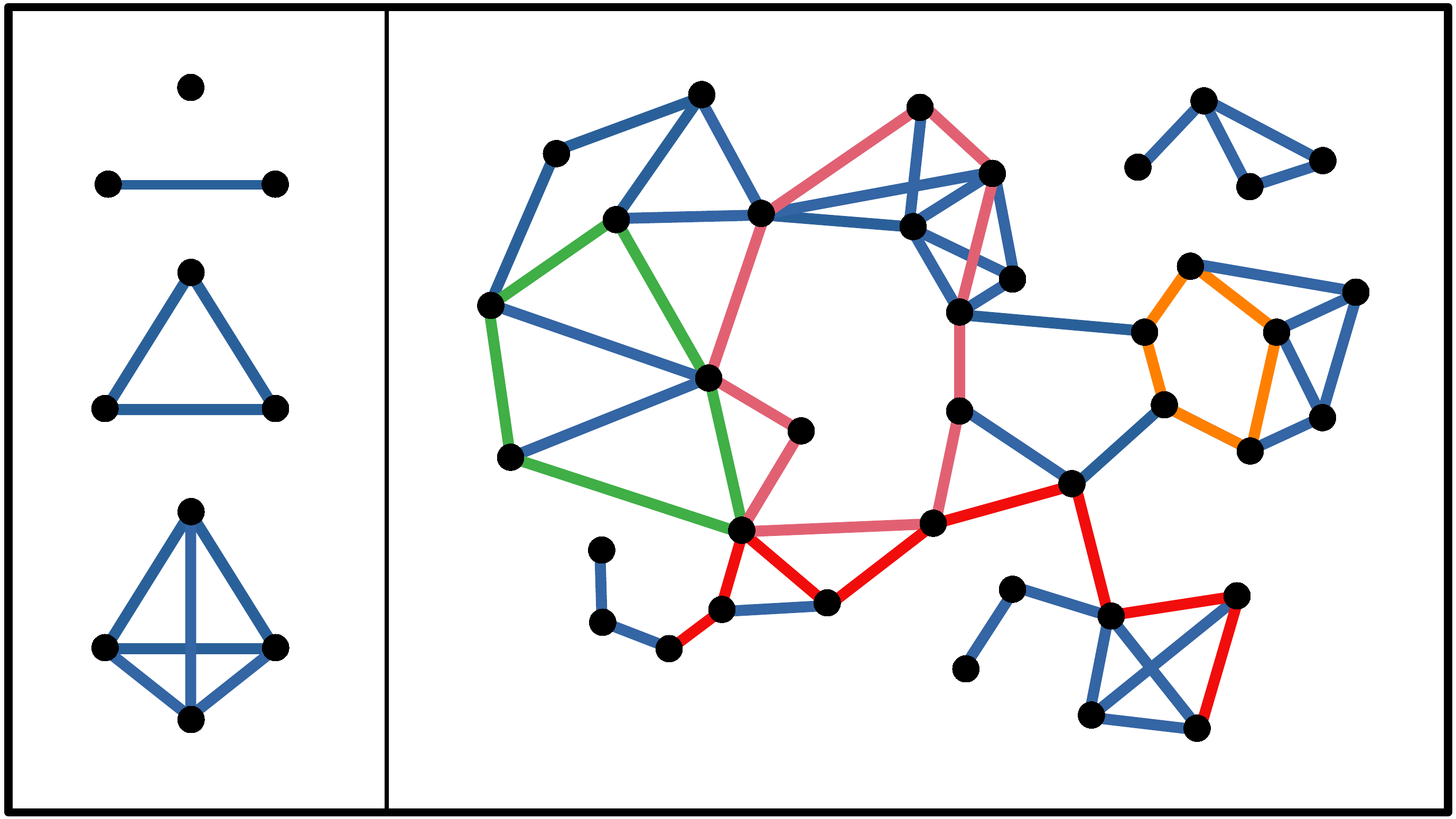}
		\caption{Left panel: low-dimensional simplicies, 0-simplex, 1-simplex, 2-simplex and 3-simplex (sorted from top to bottom). Right panel: a 3-dimensional simplicial complex containing $|\Sigma_{0}|=39$ 0-simplicies, $|\Sigma_{1}|=63$ 1-simplicies, $|\Sigma_{2}|=25$ 2-simplicies, $|\Sigma_{3}|=3$ 3-simplicies and determined by the Betti numbers $\vec{\beta}=(2,4,0,0)$. The subspaces determined by different colors illustrate example for a 1-chain (red), 1-cycle (pink), 1-boundary (green) and 1-hole (orange).}
		\label{fig:sim_com}
	\end{figure}

	\begin{figure*}
		\begin{center}
			\includegraphics[width=0.19\textwidth]{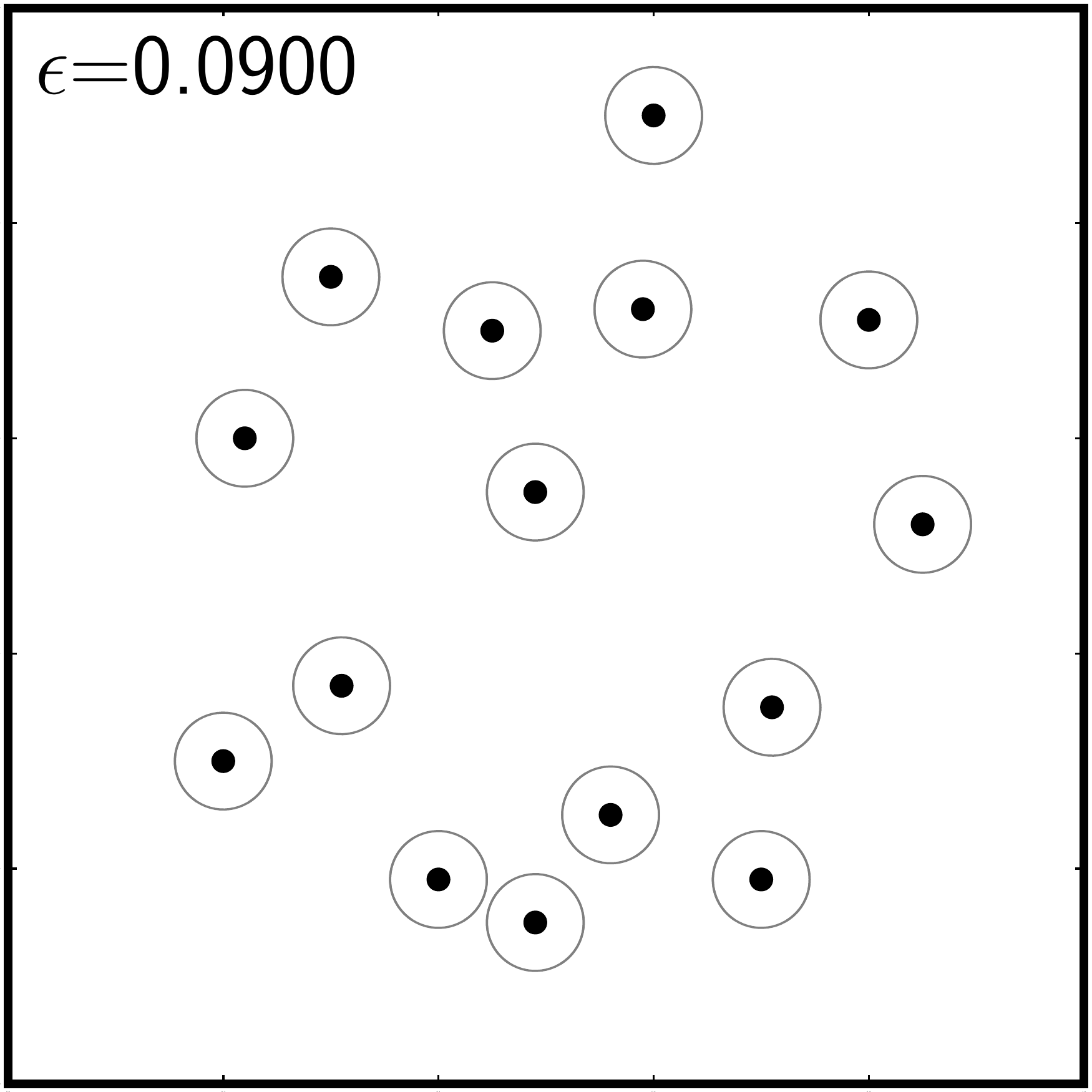}
			\includegraphics[width=0.19\textwidth]{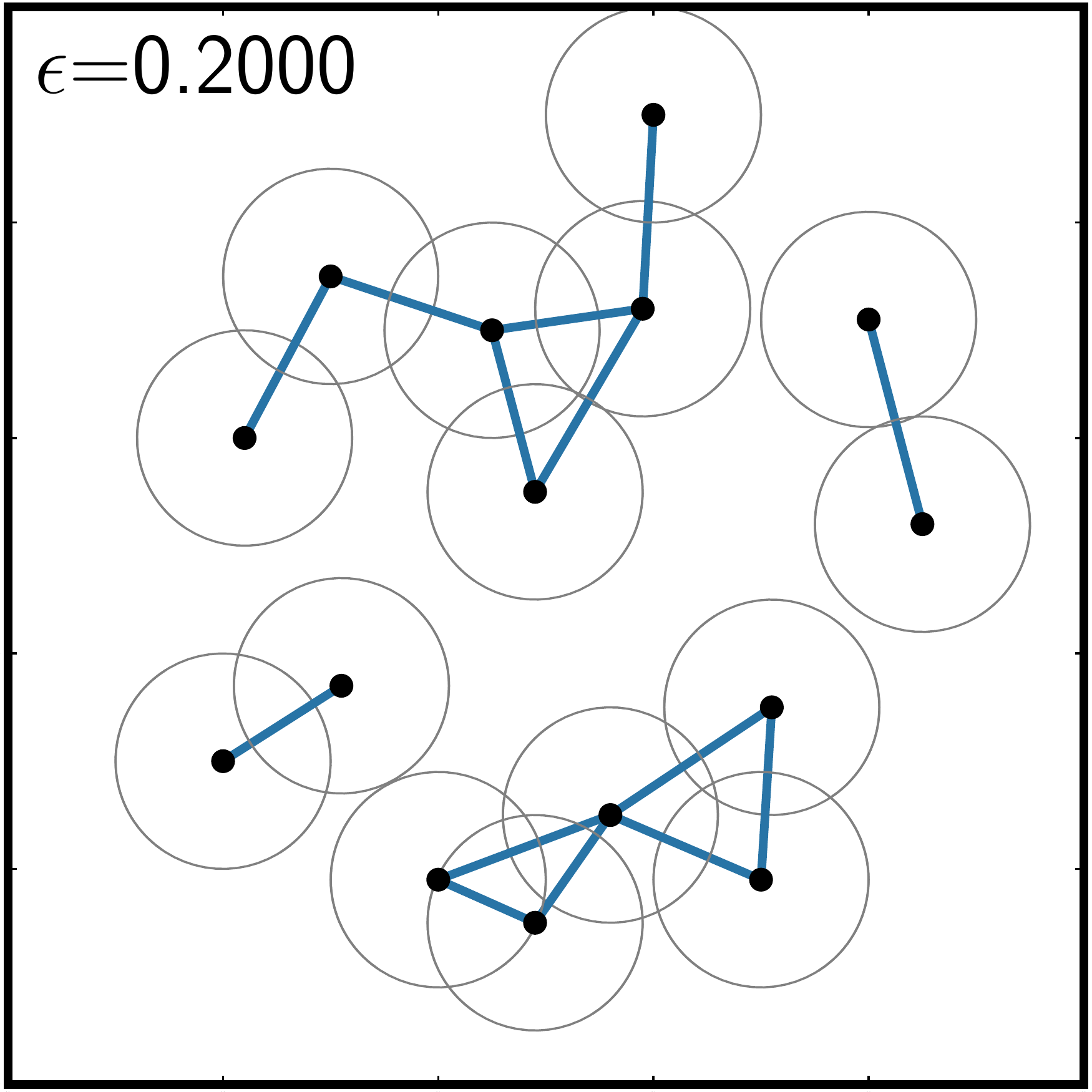}
			\includegraphics[width=0.19\textwidth]{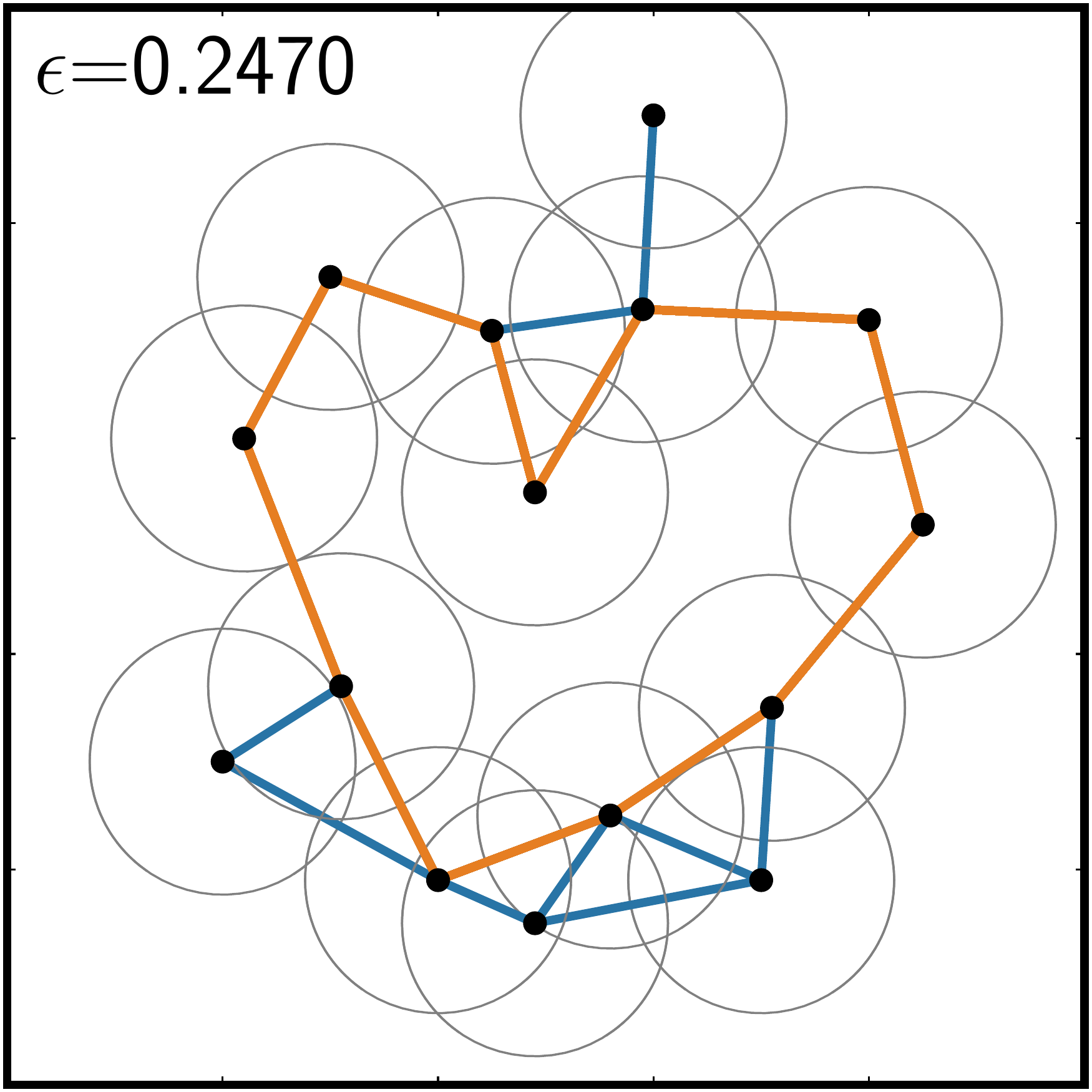}
			\includegraphics[width=0.19\textwidth]{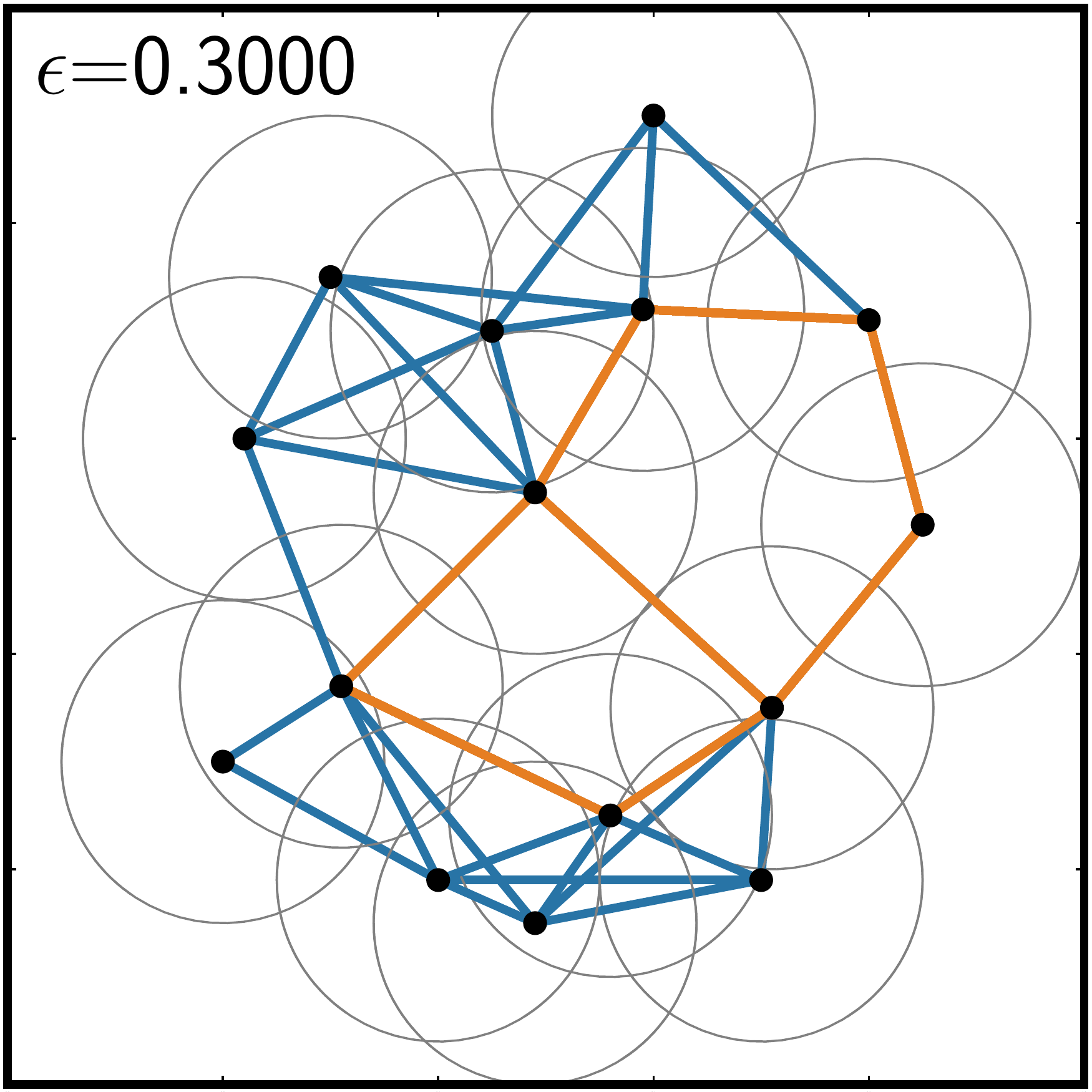}
			\includegraphics[width=0.19\textwidth]{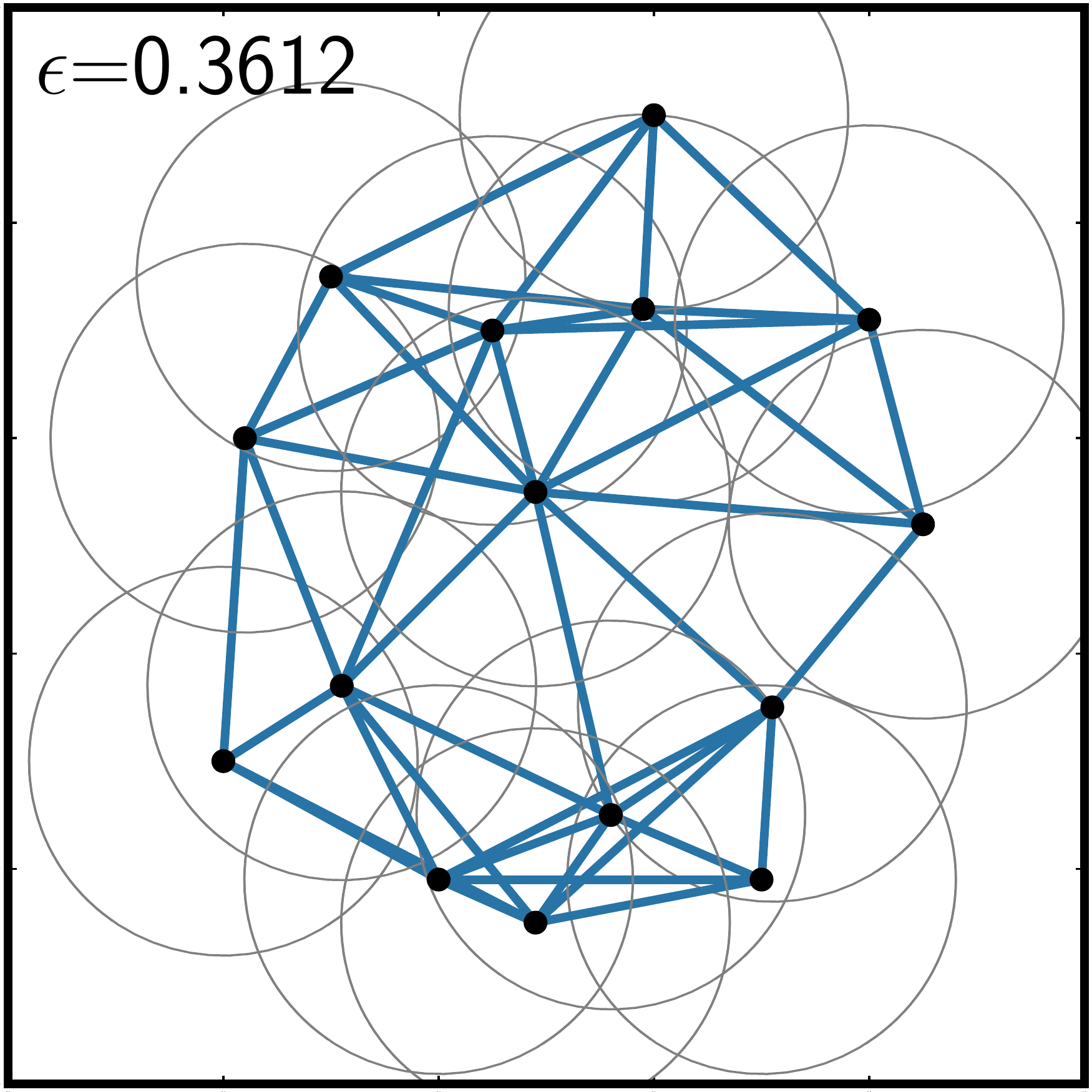}
			\includegraphics[width=0.98\textwidth]{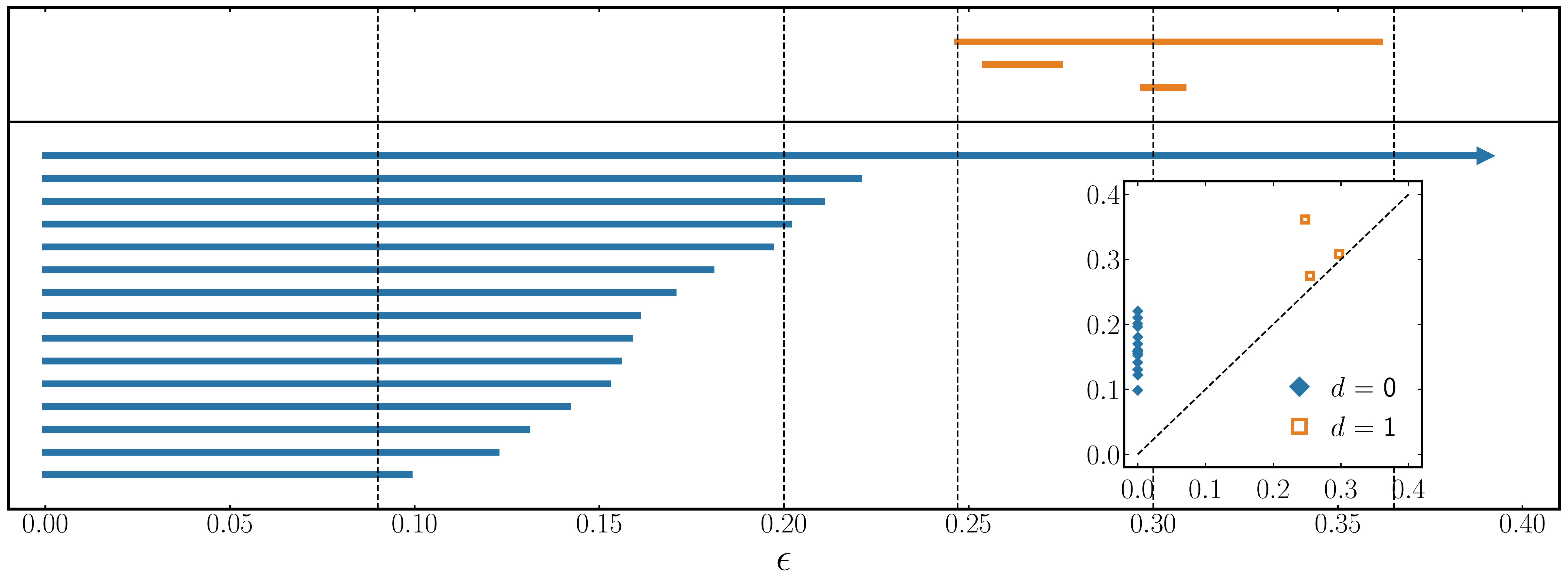}
			\caption{The persistent homology (PH) of a random PCD. This figure shows how the PH method extracts the topological features of a random 2-dimensional PCD of size $N=15$. Top panel: growing VR simplicial complex by increasing the proximity parameter $\epsilon$ (from left to right). Through the filtration, topological properties of the PCD change (connected components represented by black points and blue lines, while topological loops represented by orange lines). Bottom panel: the PB (PD in inset) associated with the filtration shows the persistence (evolution) of extracted topological features by persistence bars (pairs),  filled blue diamond and  empty orange square symboles are for the zeroth and first homology generators, respectively.}
			\label{fig:filt}
		\end{center}
	\end{figure*}

	%---------------------------------------------------
	\section{Methodology and our Pipeline}\label{method}
	
	For the sake of clarity, we will give a brief about the computational methods utilized in this paper to assess synthetic fractional Brownian motion (fBm) signal by paying attention to the mathematical preliminaries. More precisely, the PH of reconstructed phase space (embedded PCD) from fBm for different noticeable specifications and our pipeline is explained in this section.       
	
	\subsection{Time Delay Embedding}
	
	The reconstruction of phase space from a typical time series has been implemented to capture the evolution of the deterministic and chaotic dynamical systems and determine the correlations between associated quantities~\cite{packard1980geometry,myers2019persistent}. Time delay embedding (TDE) is a mathematically well-defined method to map a time series into high-dimensional Euclidean space to make such a finite-dimensional phase space~\cite{takens1981detecting}. For a given discrete time series, $x=\{x_{t}\}_{t=1}^{T}$, of length $T$, we make a set of $N$ state vectors of dimension $D$ in $D$-dimensional Euclidean space, $\mathbb{R}^{D}$, which is called $D$-dimensional PCD, for a given time-delay, $\tau$, as follows: 
	\begin{equation}
		X(x,D,\tau) = \Bigr\{ \vec x_{t} \in \mathbb{R}^{D} ~ \Bigr| ~ \vec x_{t} \equiv ( x_{t}, x_{t+\tau}, ..., x_{t+(D-1)\tau}) \Bigr\}_{t=1}^{N}
	\end{equation}
	Notice that $X(x,D=1,\tau) = x$ and the PCD size $N=T-(D-1)\tau$. According to the Takens’ theorem \cite{takens1981detecting}, we can recreate a topologically equivalent $D$-dimensional phase space from a time series by means of this method. Figure ~\ref{fig:tde} shows a synthetic fBm time series with $H=0.5$ of length $T=10^{4}$ (top panel) and associated 2-dimensional PCD of size $N=9 \times 10^{3}$ reconstructed by the TDE method for $\tau=10^{3}$ (bottom panel). As noticed before, for any value of $D$ and $\tau$, one can construct PCD from fBm, $x(H,q,T)$, and therefore, the associated homology groups are examined.

	\subsection{Persistent Homology}
	
	%\subsubsection{Homology Theory}
	
	Topology is a branch of pure mathematics dealing with the abstract objects living in high-dimensional topological spaces (e.g. the real line, sphere, torus, and more complicated spaces) to classify them in terms of their global properties (e.g. connectedness, genus, the Euler characteristics, etc.)~\cite{nakahara2003geometry}. These global features of the space usually are topologically invariant, which means that they do not change under continuous deformations and are not dependent on the way the corresponding space is made (triangulated, mathematically). By the term continuous deformation, we mean any continuous changes (elastic motions) of the space like shrinking, stretching, rotating, reflecting, etc., but not cutting or gloving. These continuous deformations are defined by the concept of \textit{homeomorphism} in algebraic topology~\cite{arnold2011intuitive}.
	
	Here we introduce the required objects and definitions in homology theory. In algebraic topology, to determine the topological properties of a space, the associated space is triangulated by mapping it into a collection of $d$-dimensional simplicies ($d$-simplicies), called simplicial complex. A $d$-simplex $\sigma_{d} = [\vec x_{0},...,\vec x_{d}]$ is a convex-hull subset of $d$-dimensional Euclidean space, $\mathbb{R}^{d}$, determined by its $(d+1)$ geometrically independent points $\{\vec x_{0},...,\vec x_{d}\}$ in $\mathbb{R}^{d}$. For instance, a point is a 0-simplex, a line segment is a 1-simplex, a triangle is a 2-simplex, a tetrahedron is a 3-simplex etc (see the left panel of Fig.~\ref{fig:sim_com}). A simplicial complex, $\mathcal{K}$, is a collection of simplices such that any subsimplex of any simplex in the simplicial complex is in the simplicial complex. Also, any pair of simplices are either disjoint or they intersect in a lower-dimensional simplex existing in the simplicial complex. The dimension of a simplicial complex is defined as the dimension of the largest simplex in it (see the right panel of Fig.~\ref{fig:sim_com})~\cite{nakahara2003geometry}. For a simplicial complex, $\mathcal{K}$, containing $|\Sigma_{d}|$ $d$-simplices, one can create a $|\Sigma_{d}|$-dimensional vector space, $\mathcal{C}_{d}(\mathcal{K})$, called $d$-chain group of the simplicial complex $\mathcal{K}$, by the basis considered as the set of all $d$-simplicies $\Sigma_{d}(\mathcal{K}) = \{\sigma_{d}^{i}\}_{i=1}^{|\Sigma_{d}|}$ of the $\mathcal{K}$ and the vectors, called $d$-chains, as fallows: 
	\begin{equation}
		c_{d} \equiv \sum_{i=1}^{|\Sigma_{d}|} a_{i} \sigma_{d}^{i} ~~ ; ~~ a_{i} ~ \in ~ \mathbb{Z}_{2} \equiv \{0,1\}
	\end{equation}
	The boundary of a $d$-simplex is the union of all its $(d-1)$-subsimplices which is obtained by applying the boundary operator $\partial_{d}$ on the simplex: 
	\begin{equation}
		\partial_{d}(\sigma_{d}) = \sum_{i=0}^{d} (-1)^{i} ~ [\vec x_{0},...,\vec x_{i-1},\vec x_{i+1},...,\vec x_{d}]
	\end{equation}
	%One can define two subspaces of the vector space $\mathcal{C}_{d}(\mathcal{K})$ based on the boundary operator. 
	The $d$-cycle group $\mathcal{Z}_{d}(\mathcal{K}) \equiv \{c_{d} \in \mathcal{C}_{d}(\mathcal{K}) ~ | ~ \partial_{d}(c_{d}) = \emptyset \}$ of $\mathcal{K}$ is defined the set of all boundaryless $d$-chains of the $\mathcal{K}$, i.e. any $d$-chain mapped to the empty space (set) is a $d$-cycle $z_{d} \in \mathcal{Z}_{d}(\mathcal{K})$. 
	Another subspace of $\mathcal{C}_{d}(\mathcal{K})$ is called $d$-boundary group $\mathcal{B}_{d}(\mathcal{K}) \equiv \{c_{d} \in \mathcal{C}_{d}(\mathcal{K}) ~ | ~ \partial_{d+1}(c_{d+1}) = c_{d} \}$ containing all $d$-chain of $\mathcal{K}$ which is the boundary of a $(d+1)$-chain in $\mathcal{C}_{d+1}(\mathcal{K})$. The elements of $\mathcal{B}_{d}(\mathcal{K})$ are called $d$-boundary. 
	Since the boundary operator satisfies the property $\partial_{d-1}(\partial_{d}(c_{d})) = \emptyset$ for any $c_{d} \in \mathcal{C}_{d}(\mathcal{K})$, we can conclude that $\mathcal{B}_{d}(\mathcal{K}) \subseteq \mathcal{Z}_{d}(\mathcal{K})$. In order to ignore $d$-cycles of the simplicial complex $\mathcal{K}$ that are also boundary, one can consider a topological equivalence relation on $\mathcal{Z}_{d}(\mathcal{K})$ such that, any pair of cycles $z_{d}^{i} , z_{d}^{j} \in \mathcal{Z}_{d}(\mathcal{K})$ are equivalent (homologous), if $z_{d}^{i} - z_{d}^{j} \in \mathcal{B}_{d}(\mathcal{K})$. This equivalence relation partitions $\mathcal{Z}_{d}(\mathcal{K})$ into a union of disjoint subsets, called $d$th homology classes. The $d$-homology
	group of the simplicial complex $\mathcal{K}$ is defined as $\mathcal{H}_{d}(\mathcal{K}) \equiv \{[z_{d} ] ~ | ~ z_{d} \in \mathcal{Z}_{d}(\mathcal{K}) \}$, where $[z_{d}]$ represents the homology class of $z_{d}$. The $d$th Betti number of simplicial complex $\mathcal{K}$, denoted by $\beta_{d}(\mathcal{K})$, as a topological invariant of $\mathcal{K}$, is the dimension of $d$-homology group of $\mathcal{K}$. Intuitively, $\beta_{d}(\mathcal{K})$ indicates the number of $d$-dimensional topological holes ($d$-holes) of the topological space triangulated by the simplicial complex $\mathcal{K}$.

	\begin{figure*}
		\includegraphics[width=\textwidth]{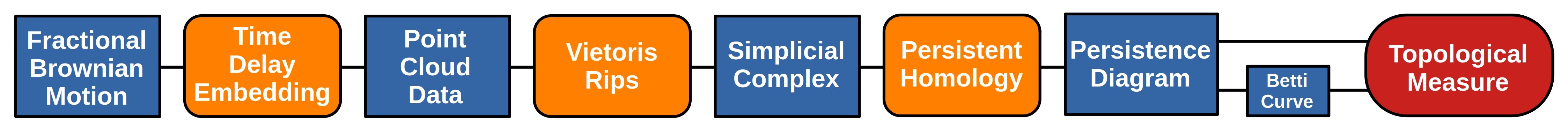}
		\caption{The proposed pipeline in this paper: A synthetic fractional Brownian motion (fBm) time series is converted to a point cloud data (PCD) by the time-delay embedding (TDE) method. The evolution of homology groups of simplicial complex mapped from the constructed point cloud data by the Vitoris-Rips (VR) method is calculated by the persistent homology (PH) technique and summarized as persistence pairs (PPs) in persistence diagrams (PDs). Finally, topological measures based on the persistence diagrams and Betti curves are computed.}
		\label{fig:pipeline}
	\end{figure*}
	
	\begin{figure*}
		\includegraphics[width=0.32\textwidth]{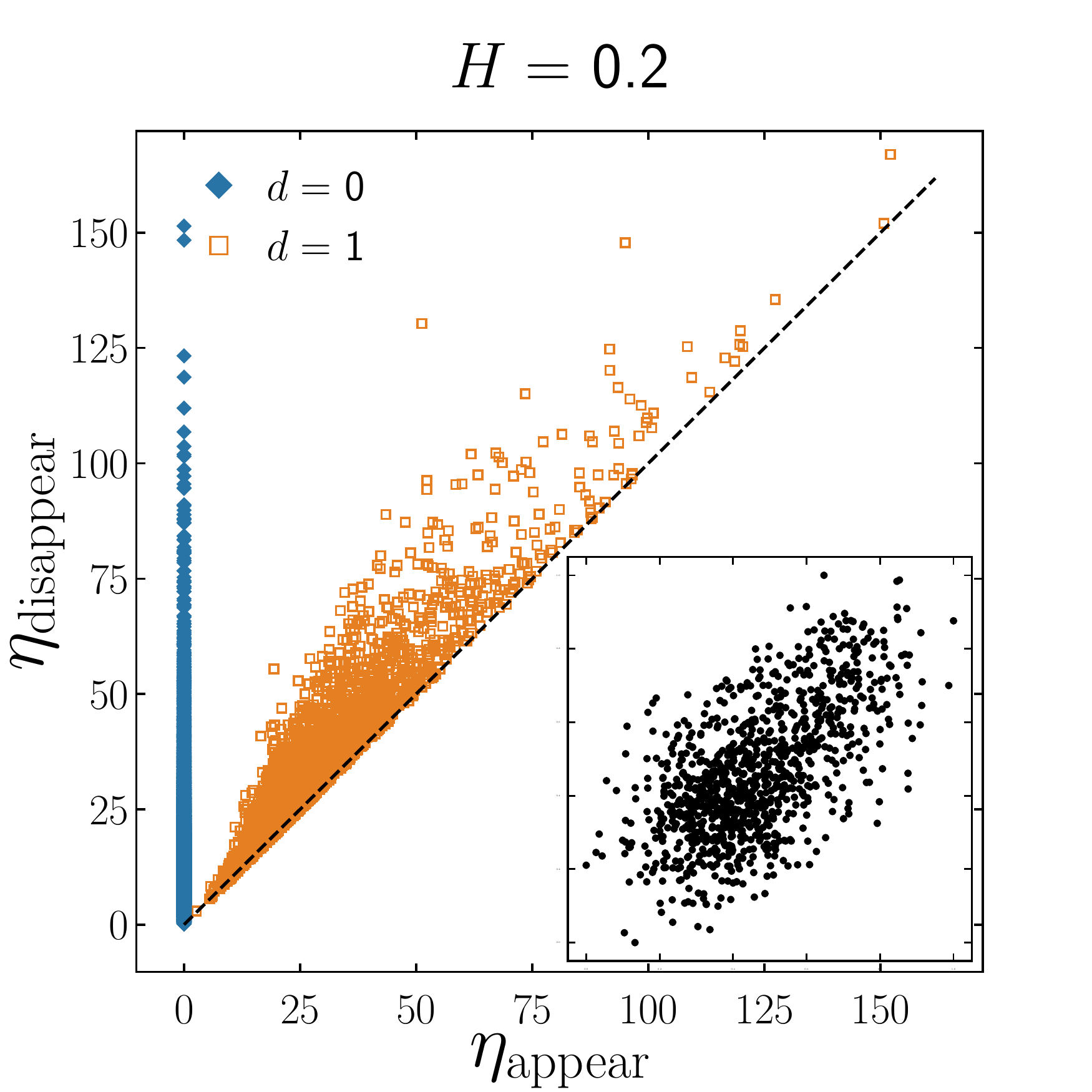}
		\includegraphics[width=0.32\textwidth]{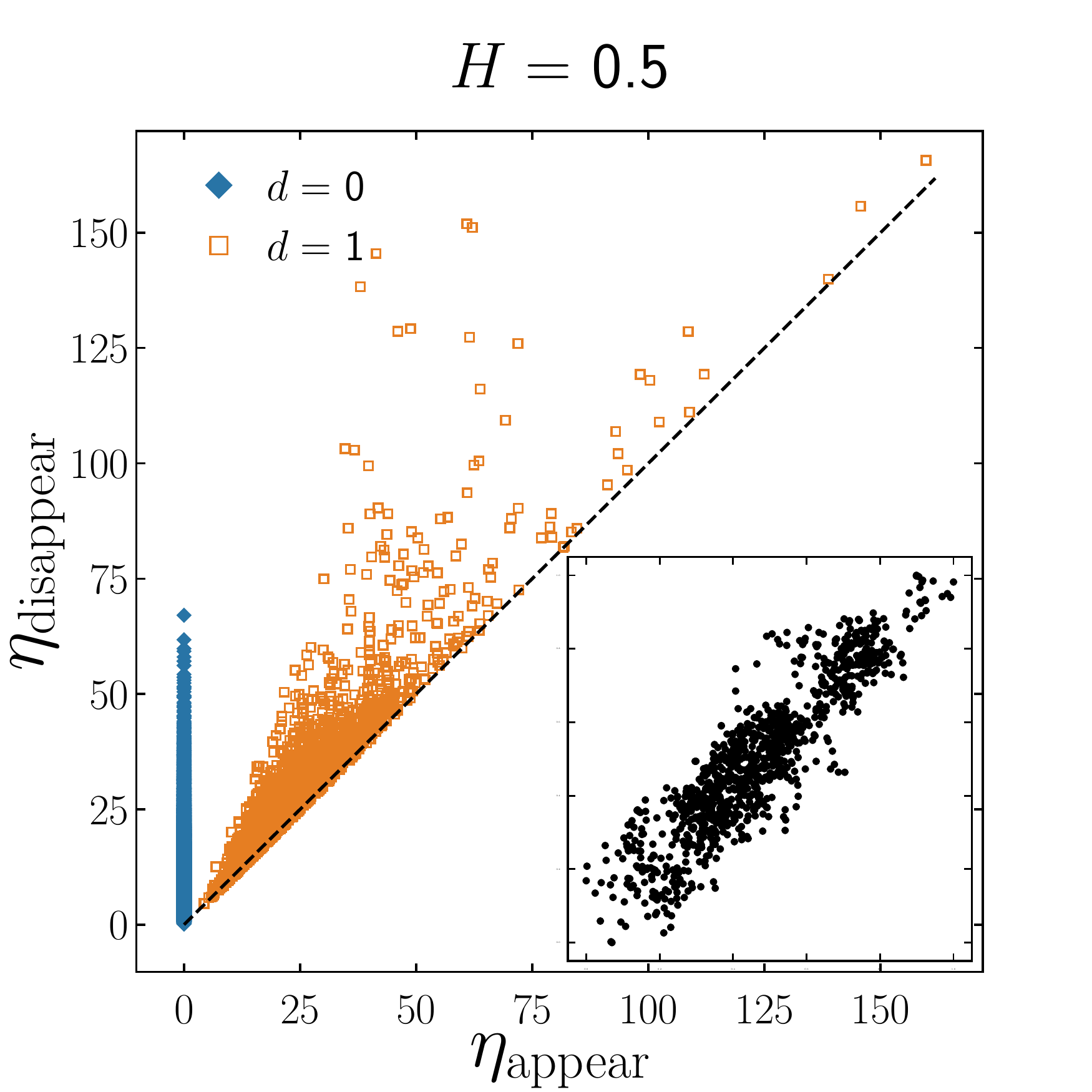}
		\includegraphics[width=0.32\textwidth]{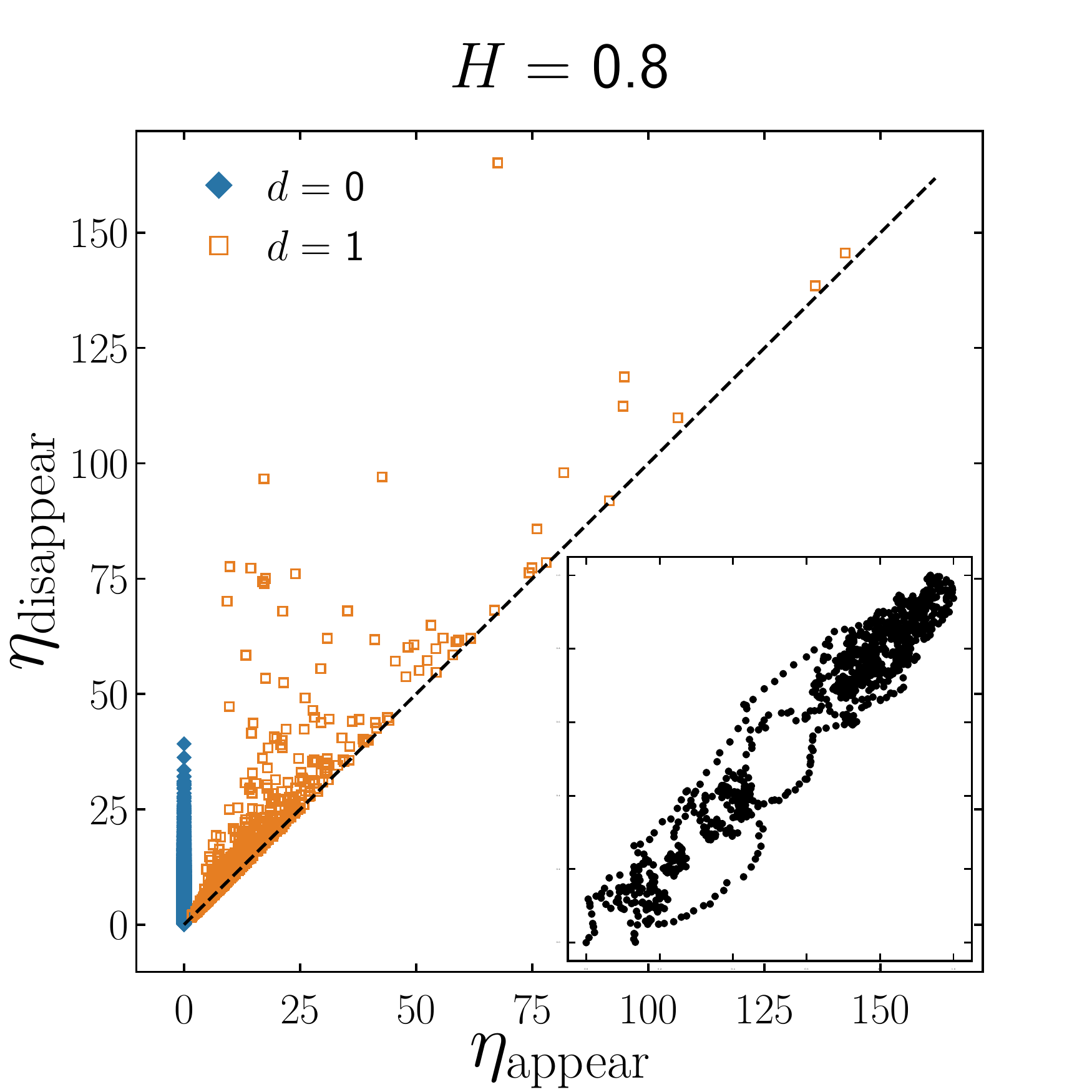}
		\includegraphics[width=0.32\textwidth]{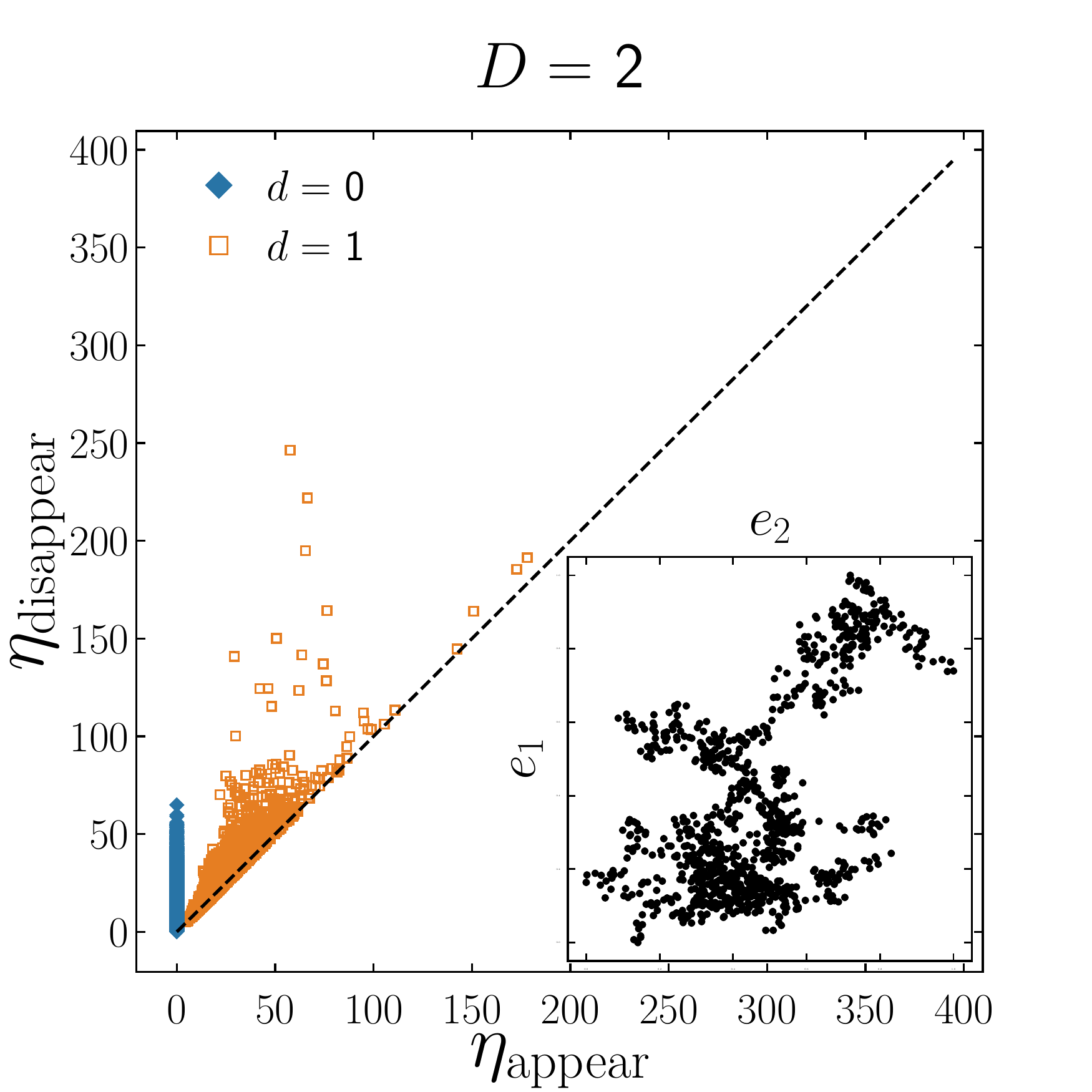}
		\includegraphics[width=0.32\textwidth]{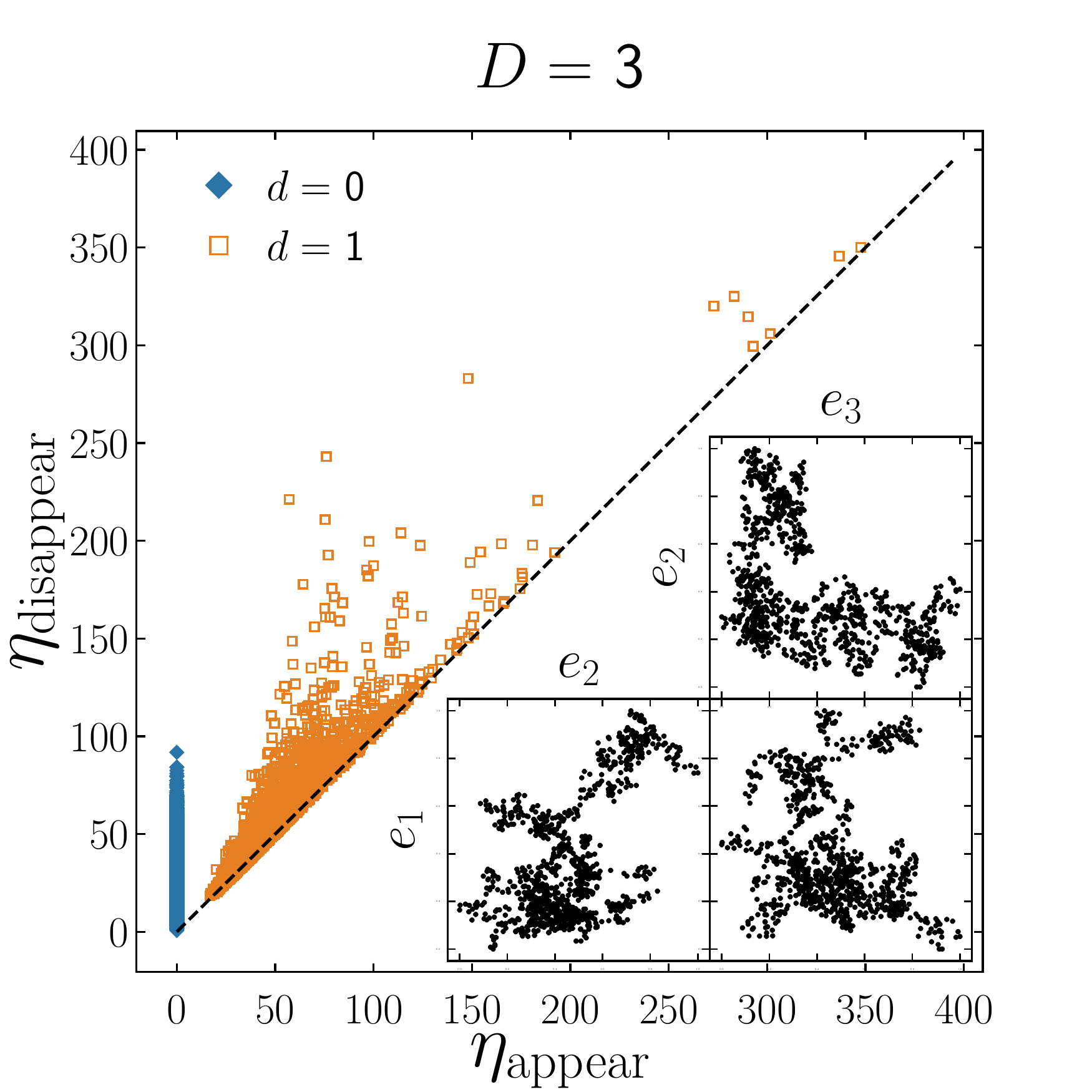}
		\includegraphics[width=0.32\textwidth]{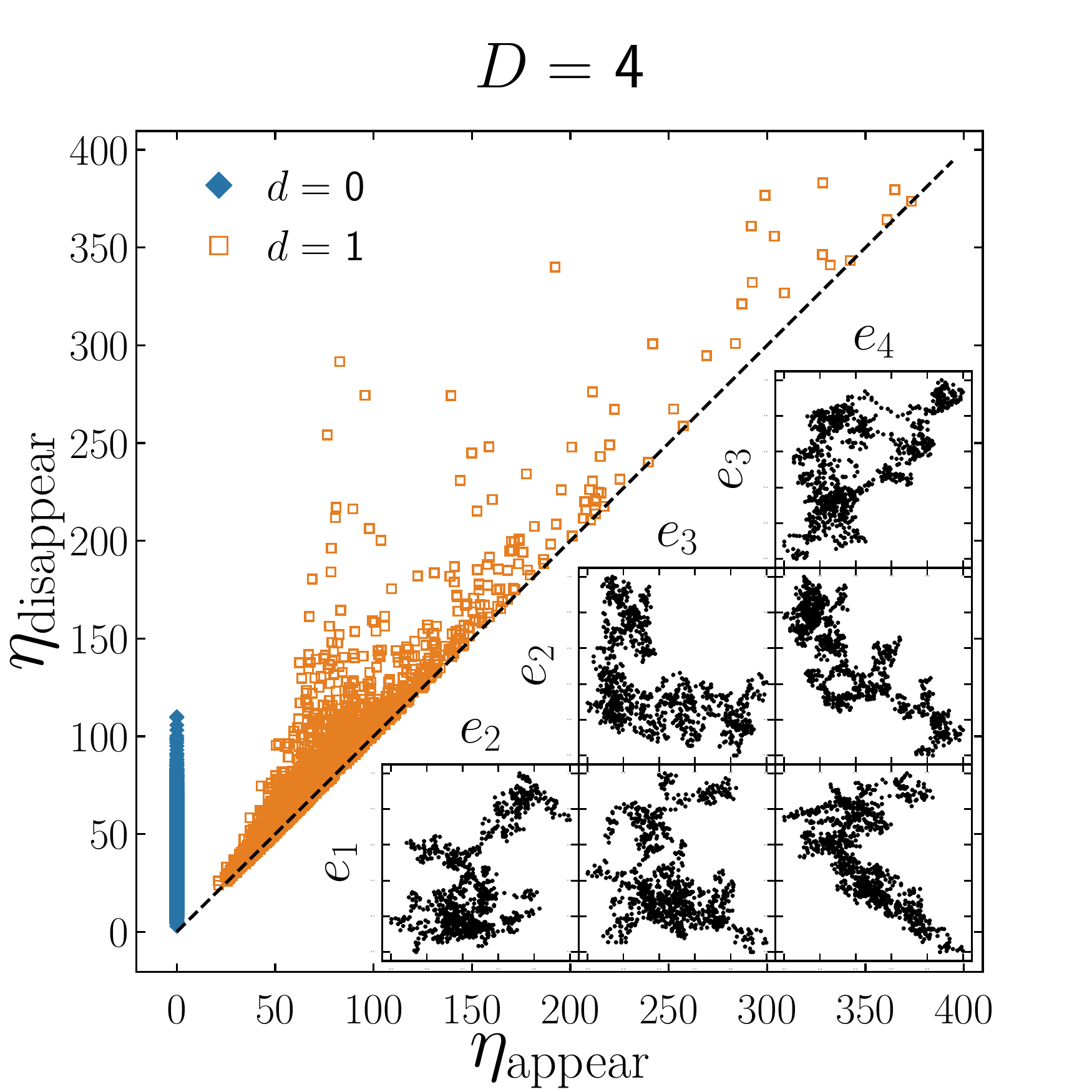}
		\caption{Upper panels: the PDs of $d$th homology groups (filled blue diamond for $d=0$ and empty orange square for $d=1$) for 2-dimensional PCD (inset plots) converted from regular fBm series with $H=0.2$ (left), $H=0.5$ (middle) and $H=0.8$ (right) by the TDE method for time-delay $\tau=100$. The population and distribution of PPs varies by $H$ which cause the $H$-dependency of topological measures. Lower panels: the $D$-dependency of $0$th (filled blue diamond) and $1$st (empty orange square) PD of regular fbm with $H=0.5$ embedded to unit $D$-cube ($D=2,3,4$ from left to right) for $\tau=1000$. The inset plots show the visualization of reconstructed PCD projected to the standard planes of the space.}
		\label{fig:H-dep_PD}
	\end{figure*}
	
	%\subsubsection{Topological Persistence}
	
	Consider a $D$-dimensional PCD $X = \{\vec x_{i} \in \mathbb{R}^{D} ~ | ~ \vec{x}_i \equiv (x_{i}^{\delta})_{\delta=1}^{D} \}_{i=1}^{N}$ of size $N$ as a finite ($N \ne \infty$) discrete subset of $D$-dimensional Euclidean space $\mathbb{R}^{D}$. The first step to study this type of dataset from topological viewpoint, i.e. to understand the topological space underlying the PCD, is finding a triangulation to tessellate the dataset~\cite{edelsbrunner2022computational}. We construct a simplicial complex $\mathcal{K}$ from the PCD $X$. The obvious triangulation of $X$ is creating a simplicial complex containing only $N$ 0-simplices, $\mathcal{K}(X) = \{\sigma_{0}^{i} \equiv [\vec x_{i}] ~ | ~ \vec x_{i} \in X \}_{i=1}^{N}$, with trivial topology ($\beta_{0}(\mathcal{K}(X)) = N, \beta_{d}(\mathcal{K}(X)) = 0$ for $d>0$). To go beyond this simple structure, one can build the simplicial complex of the PCD in larger scale, i.e. higher value of proximity of the vectors in $X$. To this end, for a fixed value of the proximity parameter $\epsilon \ge 0$, we can create a simplicial complex $\mathcal{K}(X,\epsilon)$ associated with the PCD $X$, as a collection of simplcies such that any $d$-simplex $\sigma_{d} = [\vec x_{0},...,\vec x_{d}]$ in $\mathcal{K}(X,\epsilon)$ corresponds to $(d+1)$ vectors $\{\vec x_{0},...,\vec x_{d}\}$ in $X$ and the Euclidean distance between any pair of the vectors $\{\vec x_{0},...,\vec x_{d}\}$ is less than the threshold. The constructed simplicial complex $\mathcal{K}(X,\epsilon)$ is called Vietoris-Rips (VR) simplicial complex~\cite{zomorodian2005topology}. After building the VR simplicial complex $\mathcal{K}(X,\epsilon)$ from the dataset, one can calculate the topological objects defined before. But the structural properties of $\mathcal{K}(X,\epsilon)$ are highly dependent on the chosen scale $\epsilon$. To overcome this issue, the persistent homology (PH) method, from topological data analysis (TDA), computes the evolution of the extracted topological invariants ($d$th Betti numbers) by varying the proximity parameter $\epsilon$, continuously~\cite{munch2017user}. Precisely, PH builds a growing sequence of VR simplicial complexes, called filtration, by increasing the proximity parameter and captures the structural changes of the simplicial complex in terms of Betti numbers $\beta_{d}(\epsilon)$, called Betti curve. Therefore, one can get the evolution of any $d$th homology class $[z_{d}]$ in filtered simplicial complex expressed by persistence pair (PP) $\epsilon^{[z_{d}]} = (\epsilon_{\rm appear}^{[z_{d}]},\epsilon_{\rm disappear}^{[z_{d}]})$ and summarize all PPs in a multiset $\mathcal{D}_{d} = \{ \epsilon_{i}^{[z_{d}]} \}_{i}$, known as $d$th persistence diagram (PD), where $\epsilon_{\rm appear}^{[z_{d}]}$ and $\epsilon_{\rm disappear}^{[z_{d}]}$ are the scales in which the homology class $[z_{d}]$ appears and disappears, respectively. Figure ~\ref{fig:filt} illustrates the mechanism of homology class extraction from a typical random 2-dimensional PCD of size $N=15$ by the PH method. The top row shows the filtration process in which by increasing the proximity parameter (diameter of the gray circles centered by the state vectors) the topology of the VR simplicial complex varies. The zeroth (first) homology classes are represented by black points and blue lines (orange lines). The bottom panel indicates the persistence barcode (PB) and persistence diagram (PD) (inset plot) of the zeroth (blue bars in PB and filled blue points in PD) and first (orange bars in PB and orange empty points in PD) homology classes associated with the filtration. 
	
	The number of PPs in $d$th PD is called $n_{d} = |\mathcal{D}_{d}|$ and for $d=0$, this quantity has trivial value $n_{0}=N$. The reason is correspondence between the state vectors and PPs in 0th PD.  To quantify the distribution of PPs in $d$th PD, one can calculate the Shannon entropy of the lifespans of the $d$th homology classes in $\mathcal{D}_{d}$. By the lifespan of a $d$th homology class $[z_{d}]$, we mean the positive quantity $\ell^{[z_d]} \equiv \epsilon_{\rm disappear}^{[z_{d}]} - \epsilon_{\rm appear}^{[z_{d}]}$. This measure for the entropy, so-called persistence entropy (PE), is formulated as follows:
	\begin{equation}
		E_{d} = - \displaystyle \sum_{i=1}^{n_{d}} \frac{\ell_{i}^{[z_d]}}{ \sum_{i=1}^{n_{d}} \ell_{i}^{[z_d]} } ~ \log\left( \frac{\ell_{i}^{[z_d]}}{ \sum_{i=1}^{n_{d}} \ell_{i}^{[z_d]} }\right)
	\end{equation}
	Betti curves are common visualization of PDs which indicate that how the population of homology classes (Betti numbers) evolve through the scale, $\epsilon$. In fact, the number of $d$-holes of a filtered simplicial complex for a fixed value of $\epsilon$, denoted by $\beta_{d}(\epsilon)$, is the number of PPs in $d$th PD for which $\epsilon_{(i)\rm appear}^{[z_{d}]} \le \epsilon$ and $\epsilon_{(i)\rm disappear}^{[z_{d}]} > \epsilon ~, ~ i=1,...,n_{d}$. Therefore, the calculation of $d$th Betti curve can be written as: 
	\begin{equation}
		\beta_{d}(\epsilon) = \sum_{i=1}^{n_{d}} \Theta\left(\epsilon - \epsilon_{(i)\rm appear}^{[z_{d}]}\right) ~ \Theta\left(\epsilon_{(i)\rm disappear}^{[z_{d}]} - \epsilon\right)
	\end{equation}
	It is possible to define other relevant measures based on the behavior of Betti curves e.g. $\epsilon_{d}^{\rm appear}$, $	\epsilon_{d}^{\rm diappear}$ and $\epsilon_{d}^{\rm maximize}$ which are known as critical scales for $d$th Betti curve and they read as:  
	\begin{equation}
		\epsilon_{d}^{\rm appear} \equiv {\rm max} \Bigr( \epsilon ~ \Bigr| ~ \int_{0}^{\epsilon} \beta_{d}(\epsilon') ~ {\rm d}\epsilon' = 0 \Bigr)
	\end{equation}
	and 
	\begin{equation}
		\epsilon_{d}^{\rm diappear} \equiv {\rm min} \Bigr( \epsilon ~ \Bigr| ~  \int_{\epsilon}^{+\infty} [\beta_{d}(\epsilon') - c_{d}] ~ {\rm d}\epsilon' = 0 \Bigr)
	\end{equation}
	where $c_{d=0} = 1$ and $c_{d>0} = 0$, and 
	\begin{equation}
		\epsilon_{d}^{\rm maximize} \equiv {\rm min} \Bigr( \epsilon ~ \Bigr| ~ \frac{\partial \beta_{d}(\epsilon')}{\partial \epsilon'}|_{\epsilon'=\epsilon} = 0 ~ , ~ \frac{\partial^{2} \beta_{d}(\epsilon')}{\partial \epsilon'^{2}}|_{\epsilon'=\epsilon} \le 0 \Bigr)
	\end{equation}
	According to the structure of VR simplicial complex at $\epsilon=0$, we obtain $\epsilon_{0}^{\rm appear} = \epsilon_{0}^{\rm maximize} = 0$. Furthermore, for $d=0$ the critical scale $\epsilon_{0}^{\rm disappear}$ separates the connected regime ($\epsilon \ge \epsilon_{0}^{\rm disappear}$) from the disconnected regime ($\epsilon < \epsilon_{0}^{\rm disappear}$). For $d=1$, the critical scales $\epsilon_{1}^{\rm appear}$ and $\epsilon_{1}^{\rm disappear}$ determine the loopless regimes (pre-loopless regime $\epsilon < \epsilon_{1}^{\rm appear}$ and post-loopless regime $\epsilon > \epsilon_{1}^{\rm disappear}$) and $\epsilon_{1}^{\rm maximize}$ is the scale at which the simplicial complex becomes loopful. 
	The integration of the Betti curves is the last proposed measure: 
	\begin{equation}
		B_{d} \equiv \int_{\epsilon_{d}^{\rm appear}}^{\epsilon_{d}^{\rm disappear}} \beta_{d}(\epsilon') ~ {\rm d}\epsilon'
	\end{equation}
	\begin{figure}
		\includegraphics[width=0.49\textwidth]{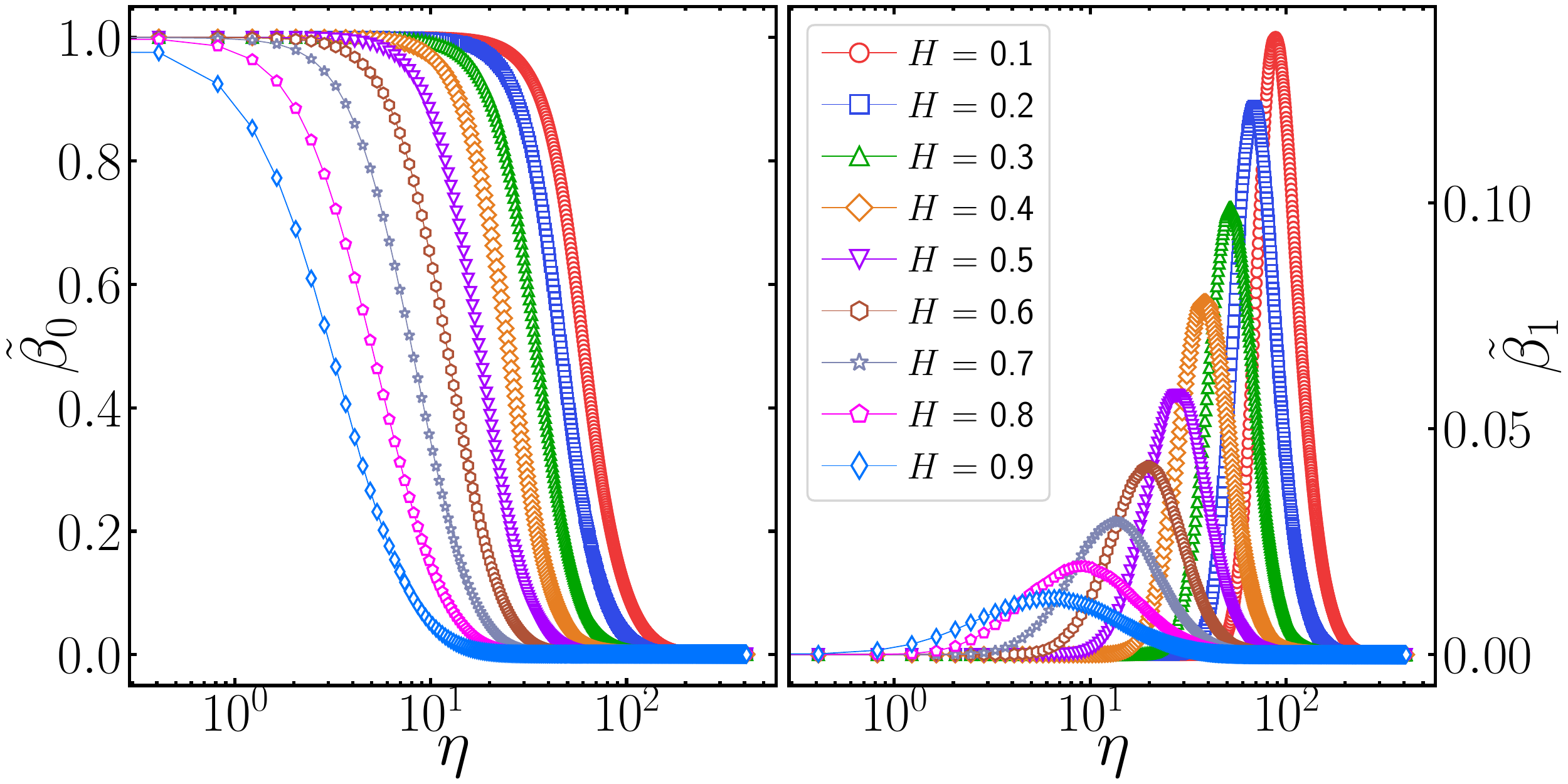}
		\includegraphics[width=0.49\textwidth]{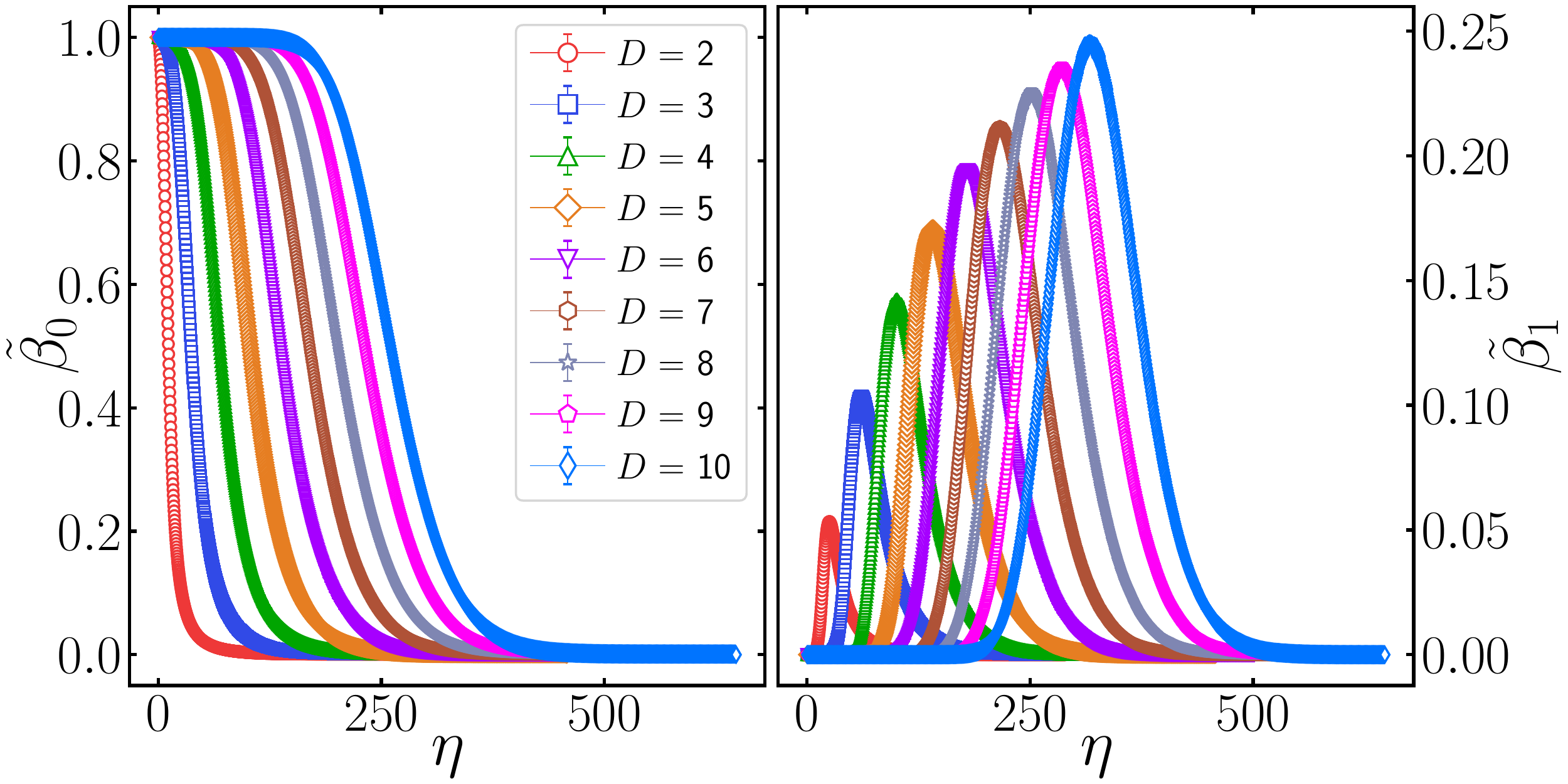}
		\caption{Upper panels: the zeroth (left) and first (lright) Betti curve of reconstructed 4-dimensional phase space for time-delay $\tau=1$ from regular fBms of various Hurst exponent. Lower panels: the Betti curves of reconstructed phase space of various dimensions by time-delay $\tau=10$ from regular fBm series with $H=0.1$. These curves indicate the strong $D$-dependency of topological properties of the embedded PCD. 
			%The topological distribution of the embedded PCD strongly depends on the Hurst exponent of the associated fBm. The PCD corresponding to the correlated series become path-connected at low scales and also the critical distances related to the first homology group is less for the signals in correlated regime.
		}
		\label{fig:betti_curve_H-dep}
	\end{figure}
	
	Now we are interested in analyzing the effect of the parameters $D$, $\tau$ and $q$ on the statistics of the $d$th homology group ($d=0,1$) and the $H$-dependency of the mentioned nontrivial topological measures, namely  $\epsilon_{0}^{\rm disappear}$, $B_{0}$, $E_{0}$, $\epsilon_{1}^{\rm appear}$, $\epsilon_{1}^{\rm disappear}$, $\epsilon_{1}^{\rm maximize}$, $B_{1}$, $E_{1}$ and $n_{1}$ and to  look for which ones are more appropriate for estimating the Hurst exponent of the fBm series in the presence of irregularity in details in next section. \\

	\begin{figure*}
		\includegraphics[width=0.95\textwidth]{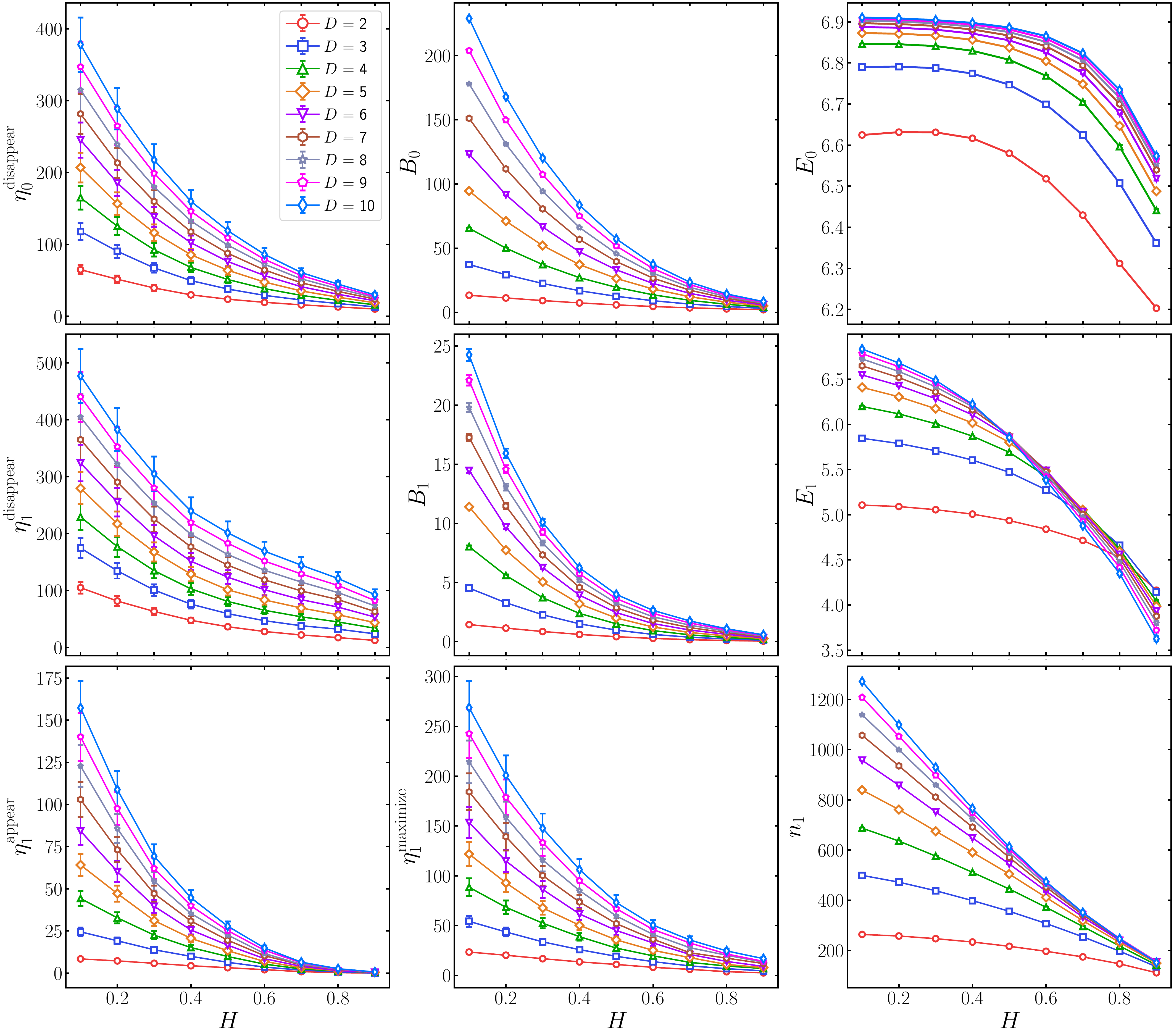}
		\caption{The $D$-dependency of topological measures introduced in the text. Here we set the parameters as  $\tau=1$ and $q=0$.}
		\label{fig:features_D-dep}
	\end{figure*} 
	\subsection{Pipeline} \label{pipeline}
	Our proposed pipeline includes the following steps (see Fig.~\ref{fig:pipeline}):
	
	(i) The fBm series $x(H,q,T)$ of distinct Hurst exponent $H$ with some irregularity $q$ with size $T$ is simulated and embedded to a high-dimensional Euclidean space by TDE method for embedding dimension $D$ and time-delay $\tau$ to construct a phase space $X(H,q,N,D,\tau)$ of size $N=T-(D-1)\tau$. 
	
	(ii) The reconstructed phase space is mapped to simplicial complex $\mathcal{K}(H,q,N,D,\tau,\epsilon)$ by the Vietoris-Rips (VR) method for various scales $\epsilon$.
	
	(iii) PH method is applied to extract the statistics and evolution of $d$-dimensional topological holes ($d$-holes) of the scale-dependent VR simplicial complex, then the associated generators are visualized as PPs in $d$th PD $\mathcal{D}_{d}(H,q,N,D,\tau)$ for zeroth and first homology groups.
	
	(iv) Some topological measures are directly computed from the PDs and some others can be defined by the criticality of the behavior of $d$th Betti curve $\beta_{d}(H,q,N,D,\tau,\epsilon)$ calculated from PDs. 
	
	It is worth noticing that, we define the normalized quantities $\tilde\beta_{d} \equiv \beta_{d}/N$ and $\eta \equiv N\epsilon$ and use them to calculate the introduced topological measurements instead of the quantities $\beta_{d}$ and $\epsilon$ in the rest of this paper.

	\begin{figure}
		\includegraphics[width=0.35\textwidth]{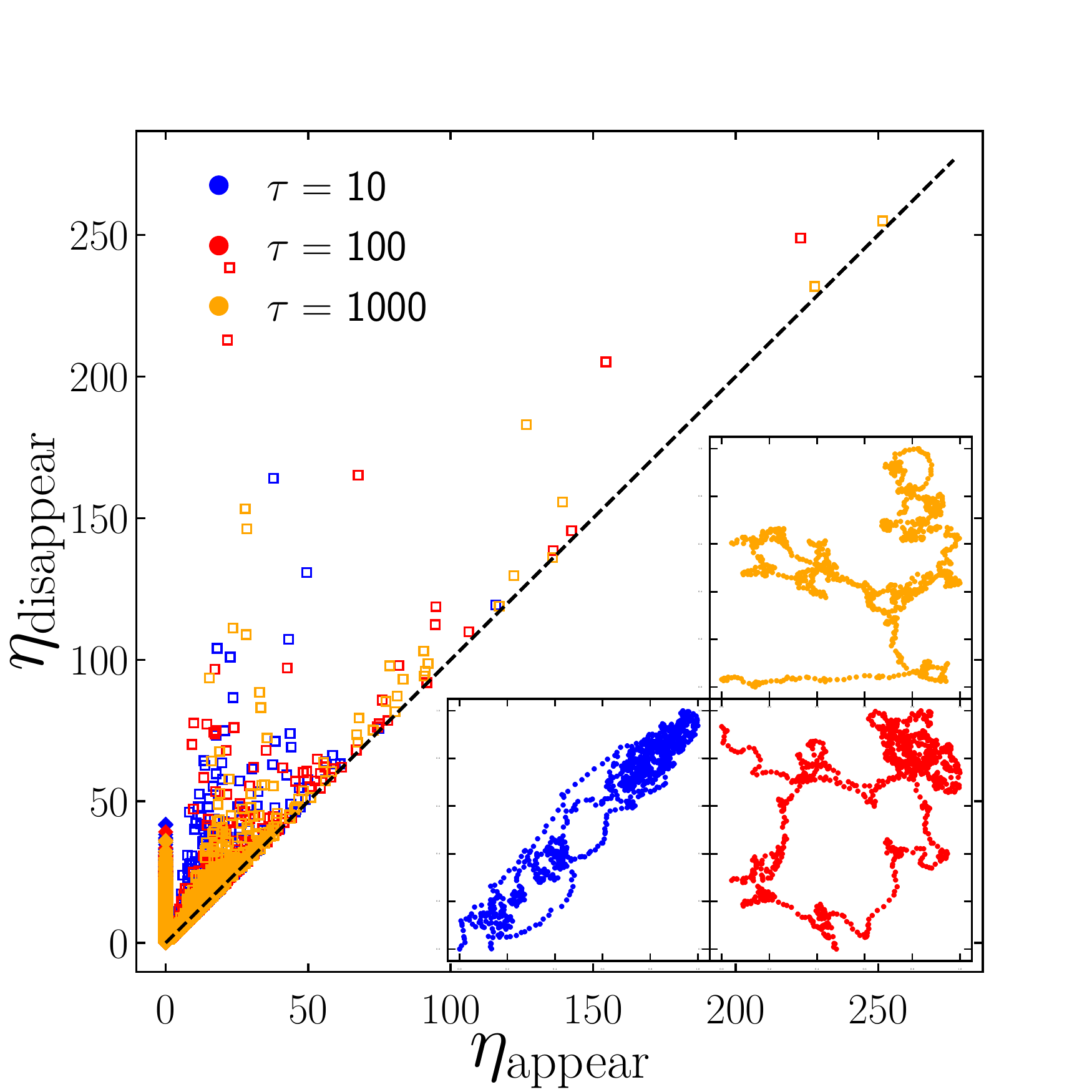}
		\includegraphics[width=0.35\textwidth]{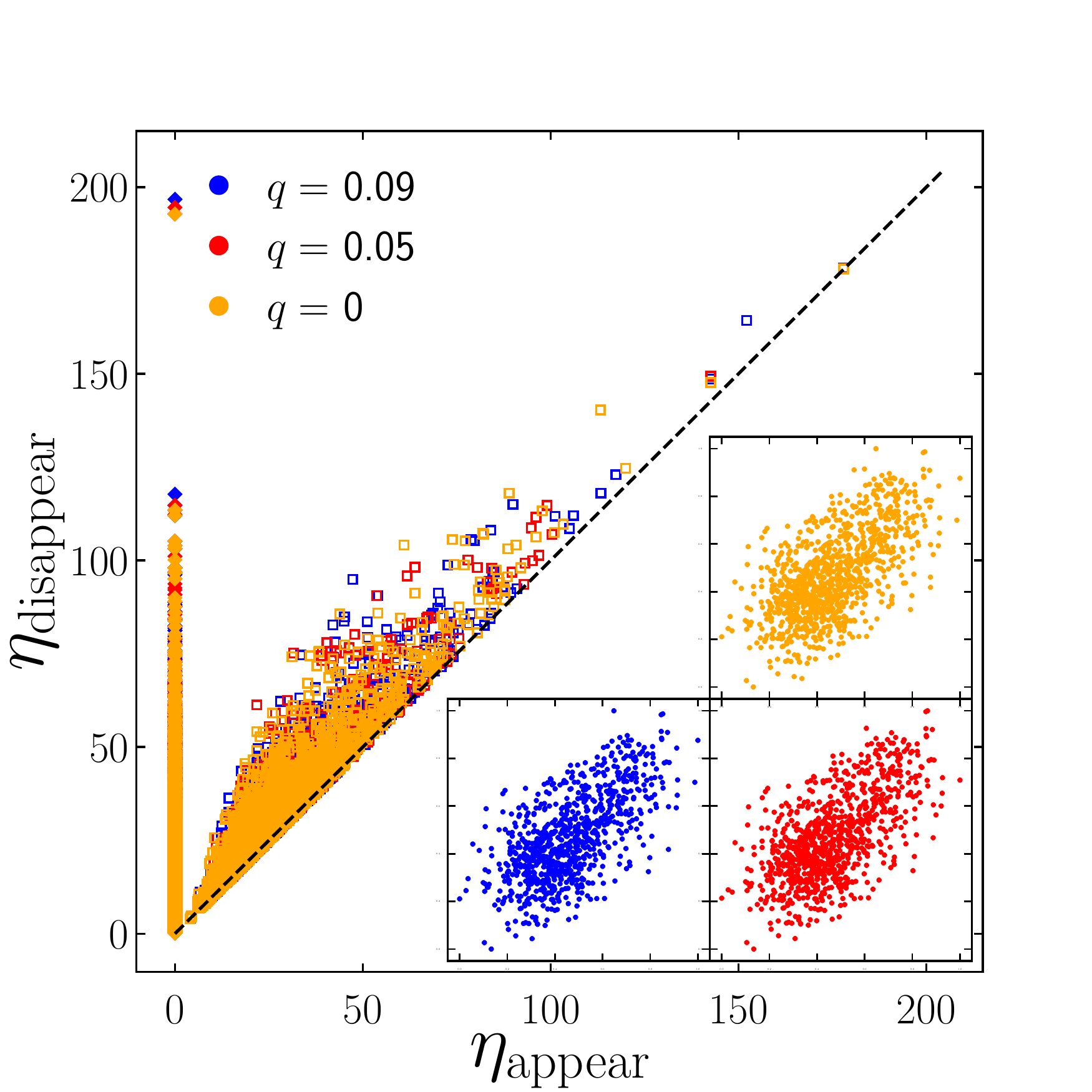}
		\caption{Upper panel: the 2-dimensional PCD of regular fBm with $H=0.9$ for time-delay $\tau=10$ (blue), $\tau=100$ (red) and $\tau=1000$ (orange) and associated PDs (inset plots). Lower panel: the PDs of 2-dimensional PCDs (insets) for $q=0$ (orange), $q=0.05$ (red) and $q=0.09$ (blue). Here we took $H=0.1$. In these plots, the associated zeroth and first PDs are represented by filled blue diamond and empty orange square symboles, respectively.}
		\label{fig:tau-dep_PD}
	\end{figure}
	
	\begin{figure}
		\includegraphics[width=0.49\textwidth]{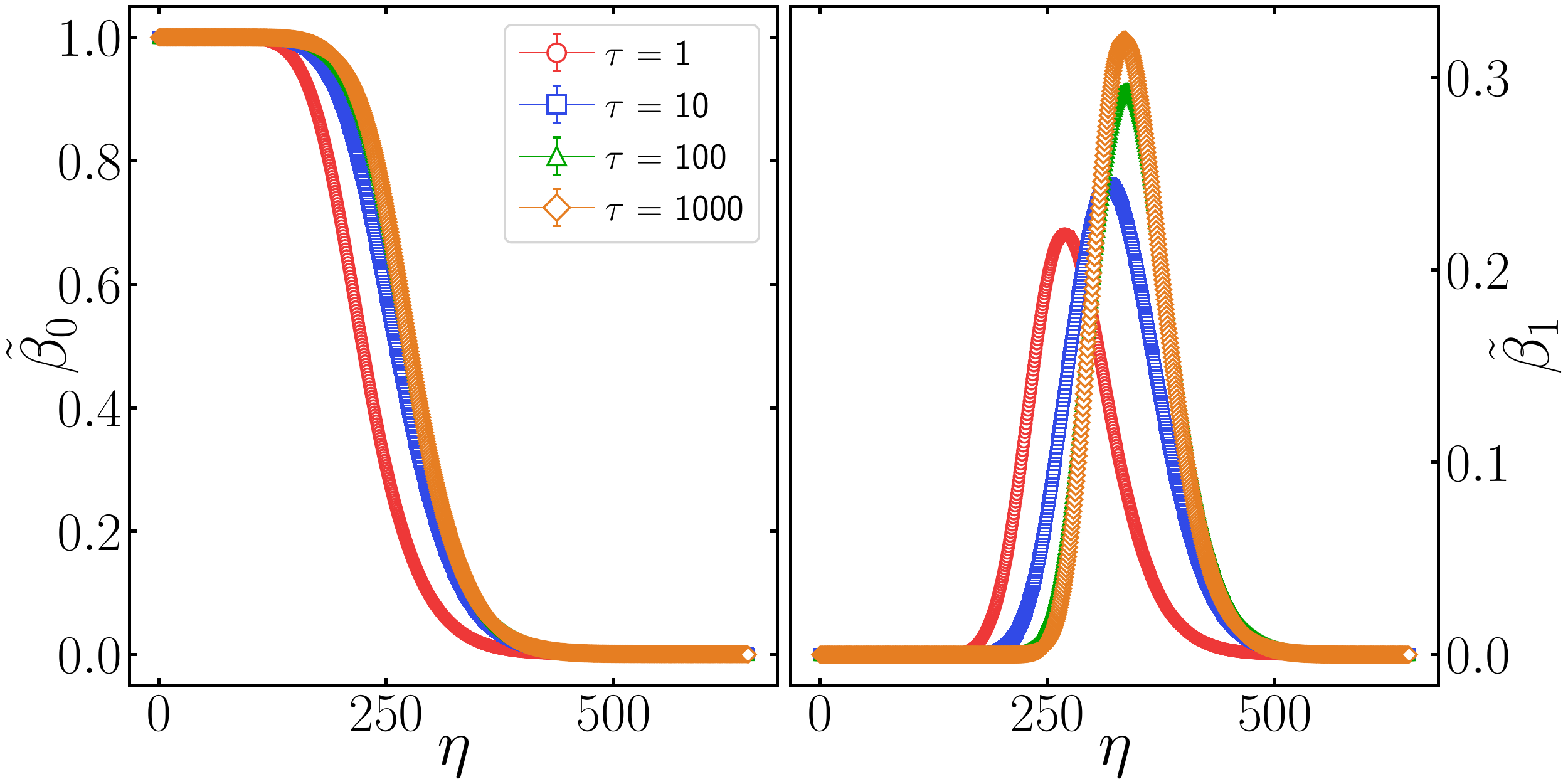}
		\includegraphics[width=0.49\textwidth]{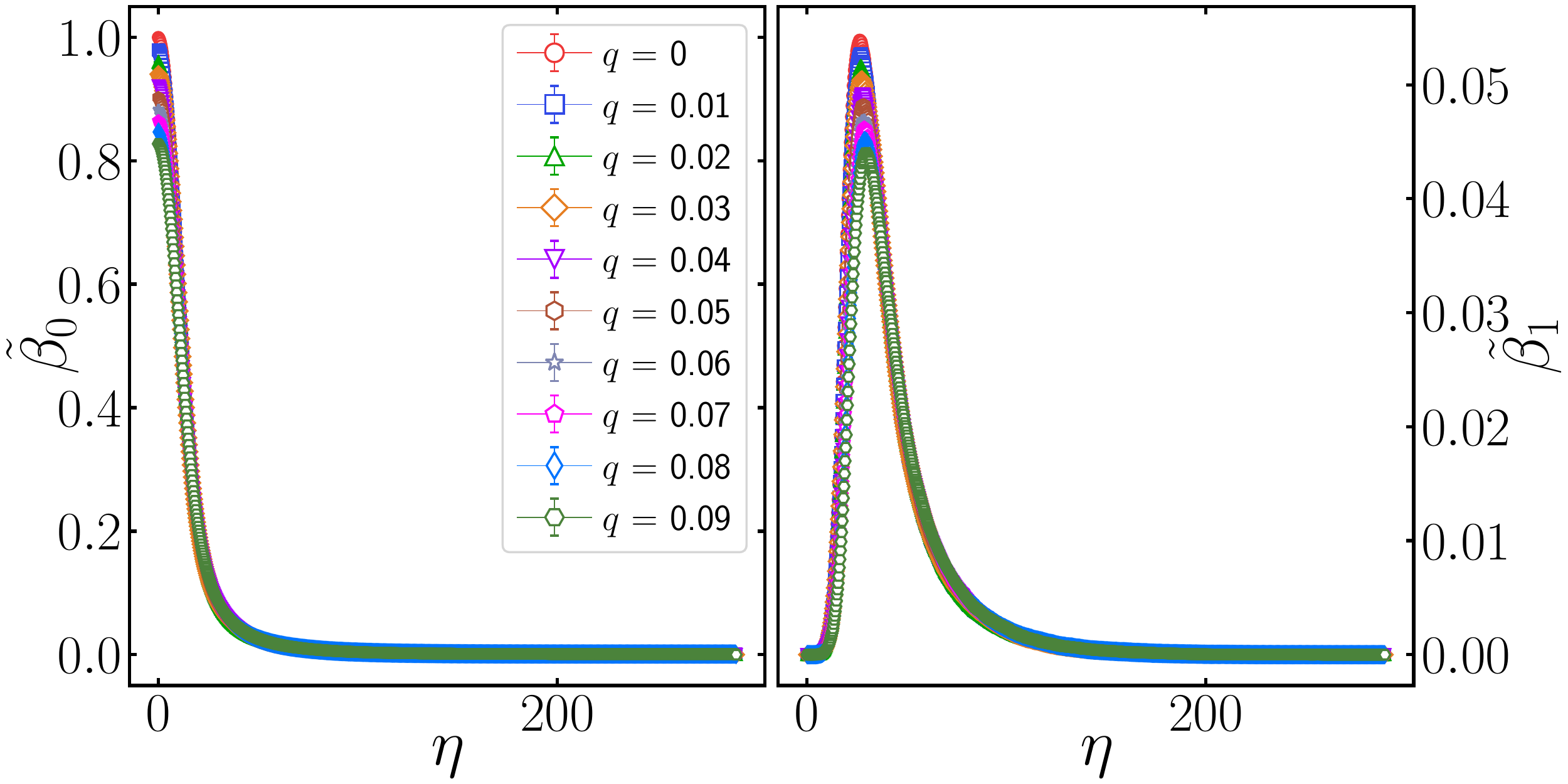}
		\includegraphics[width=0.49\textwidth]{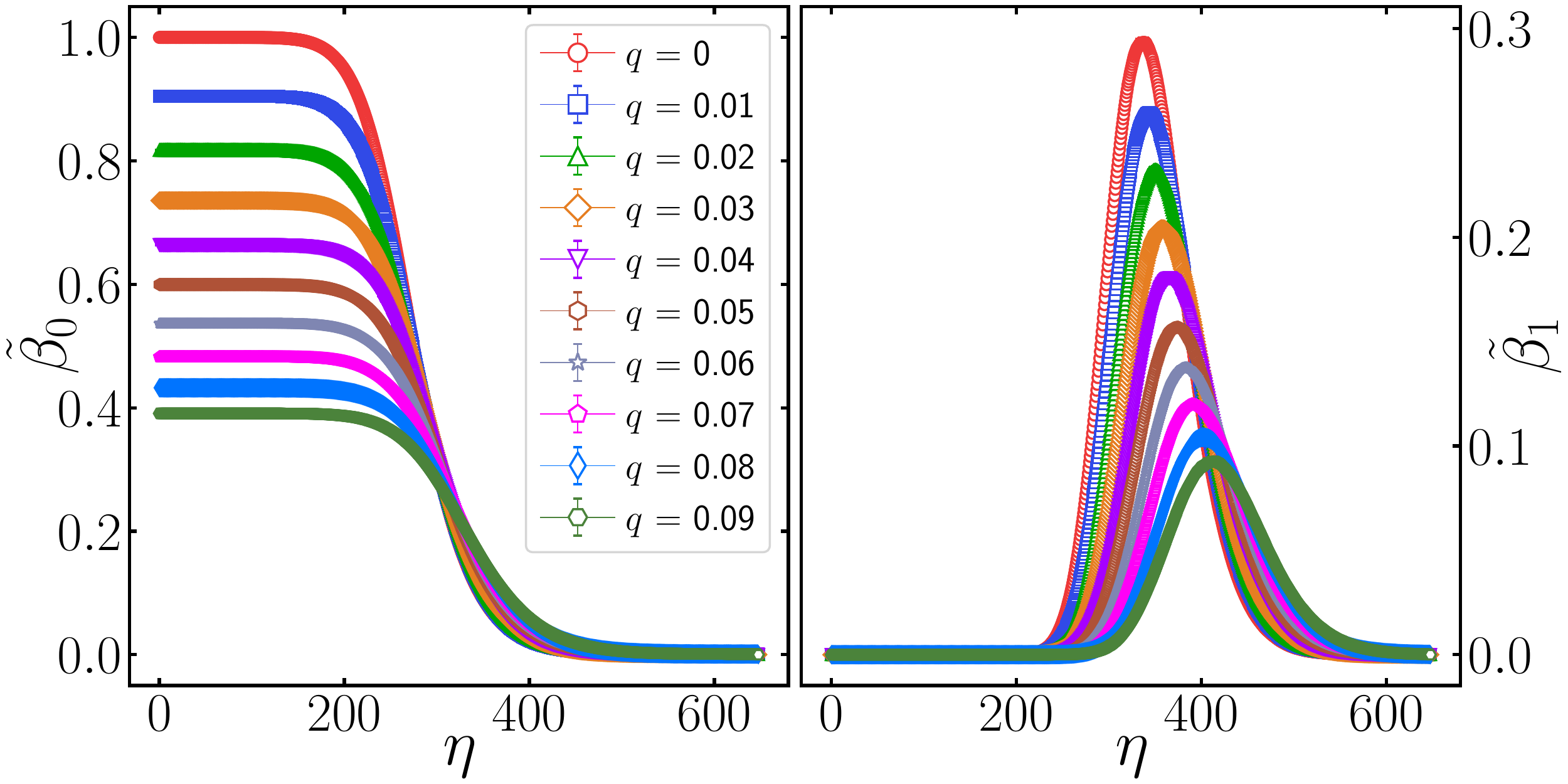}
		\caption{Upper panels: the evolution of $0$th (left) and $1$st (right) Betti number of reconstructed 10-dimensional phase space from regular fBm of $H=0.1$ for different time-delay parameter. Middel panels: the number of connected components (left) and topological loops (right) for $D=2$ and $H=0.1$ as a function of threshold for various value of irregularity $q$. The lower panel is the same as middle panel except for $D=10$.}
		\label{fig:betti_curve_tau-dep}
	\end{figure}
	
	\section{Results}\label{results}
	
	According to our pipeline, at first, we generate synthetic fBm for any given Hurst exponent based on  the Holmgren-Riemann-Liouville fractional integral \cite{mandelbrot1968fractional,kahane1993some,reed1995spectral}
	\begin{eqnarray}
		x(H,t)=\frac{1}{\Gamma\left(H+\frac{1}{2}\right)}\int_0^t(t-s)^{H-\frac{1}{2}}{\rm d}x(s)
	\end{eqnarray}
	here $\Gamma$ is the Gamma function and ${\rm d}x(s)=x(s+{\rm d}s)-x(s)$ is the increment of $(1+1)$-dimensional Brownian motion (see \cite{masoomy2021persistent} for more details). Then  we convert and rescale our mock fBm time series of various Hurst exponent $H \in (0,1)$ with $\Delta H=0.1$ to the unit $D$-dimensional cube ($D$-cube) $[0,1]^{D}$, a subspace of $D$-dimensional Euclidean space, from $D=2$ to $D=10$ ($\Delta D=1$), by using TDE method for $\tau=1,10,100,1000$ and irregularity value, defined as the fraction of missing data points in time series to the time series length $q \equiv T_{\rm missing} / T$, from $q=0$ (regular) to $q=0.09$ by step $\Delta q = 0.01$. Notice that we have valued the length of time series $T$, such that the PCD size $N=T-(D-1)\tau$ be fixed $N=2^{10}$ and the missing data points are selected uniformly randomly. In PH part, the proximity parameter varies continuously from $\epsilon_{\rm min}=0$ to $\epsilon_{\rm max}=0.2\sqrt{D}$, where $\sqrt{D}$ is the maximum possible distance between any two points in unit $D$-cube (diameter). All proposed topological measurements in this work are computed and averaged over $10^{3}$ realizations. For computational part, we utilize \textit{"Ripser"} Python package~\cite{ctralie2018ripser}.

	Now, we are going to evaluate topological measures for reconstructed phase space (embedded PCD) from generated fBm. We are interested in examining the effect of relevant parameters, namely $(H,q)$ (so-called intrinsic parameters) and $(D,\tau)$ (algorithmic parameters) on $d$th homology group ($d=0,1$). The upper panels of Fig.~\ref{fig:H-dep_PD} show the $d$th PD ($d=0$ filled blue diamond and $d=1$ empty orange square) for Hurst exponents of $H=0.2, 0.5, 0.8$ (from left to right) for fixed parameters $D=2$, $\tau=100$ and $q=0$. The insets are the visualization of the reconstructed PCD. The $d$th Betti curves of various $H$ are also shown in the upper panels of  Fig.~\ref{fig:betti_curve_H-dep}. In this part, the parameters are $D=4$, $\tau=1$ and $q=0$. Mentioned plots reveal how the topological distribution of state vectors in $D$-cube and consequently the corresponding PDs and Betti curves change for various values of Hurst exponent of the fBm series. The amount of memory encoded in the correlation quantity of the fBm signal impacts the pattern of reconstructed PCD, such that the higher value of $H$ corresponding to a higher value of correlation in fBm leads to the continuous trajectory (slow changes) than the lower value of correlation. In fact, by decreasing $H$, the $N$ state vectors become randomly distributed in $D$-cube.   Almost all the proposed topological measures decrease by increasing the Hurst exponent for the various value of parameters $D$, $\tau$, and $q$. 

	\begin{figure*}
		\includegraphics[width=0.95\textwidth]{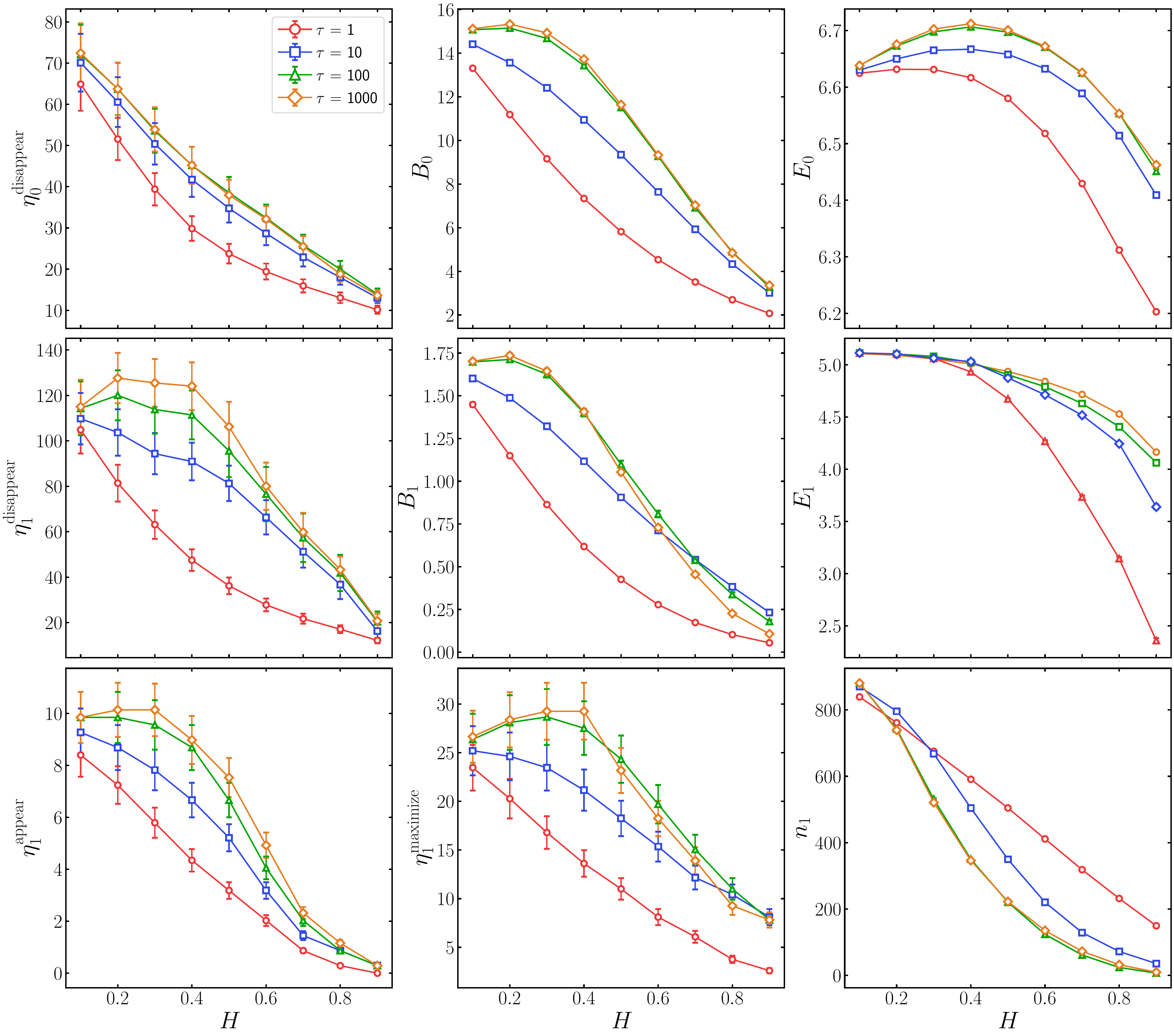}
		\caption{The time-delay parameter $\tau$ impacts on the $H$-dependency of topological measures.}
		\label{fig:features_tau-dep}
	\end{figure*}
	
	\begin{figure*}
		\includegraphics[width=0.95\textwidth]{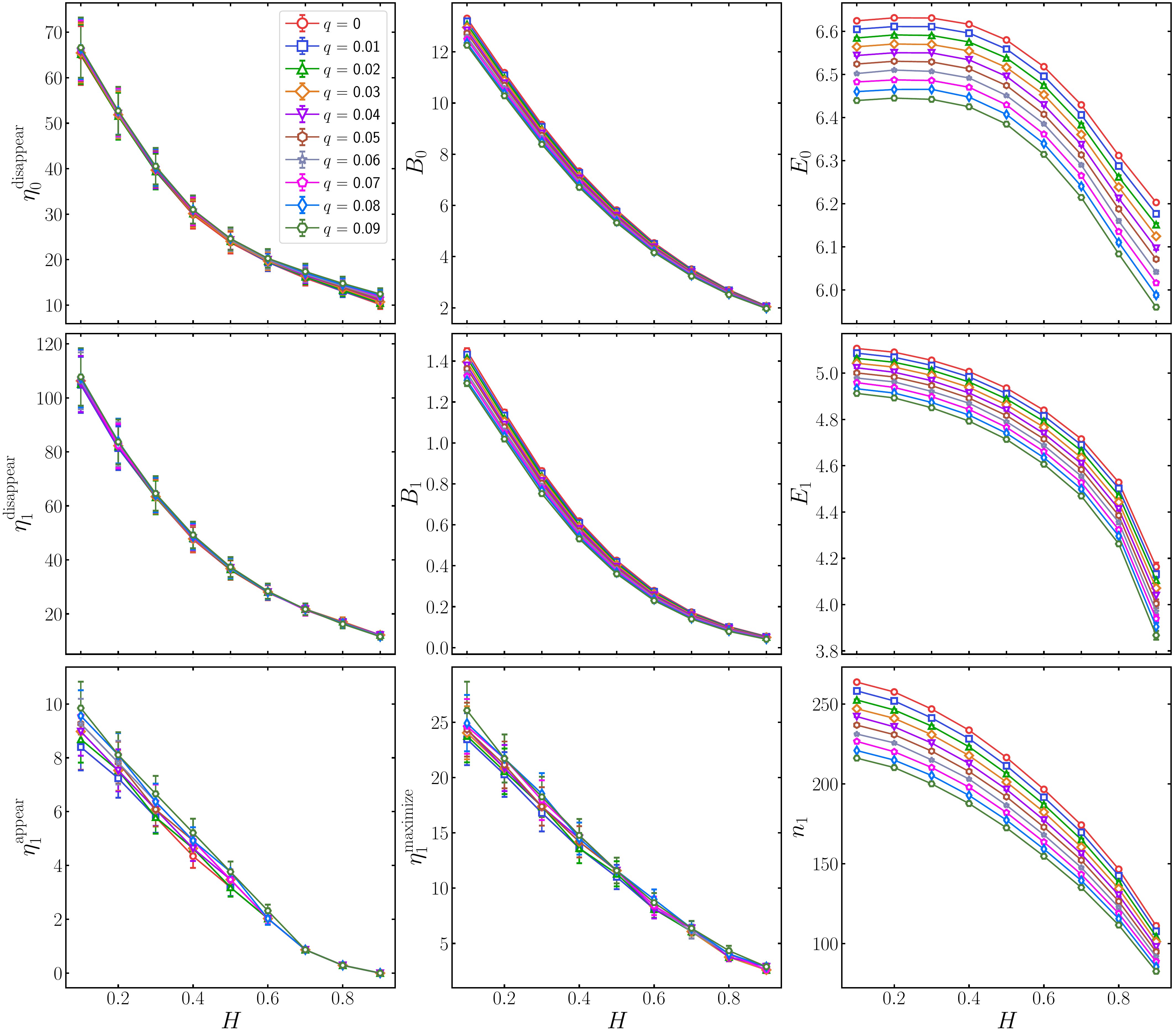}
		\caption{The proposed measures as a function of the Hurst exponent for various amounts of irregularity controlled by $q$.}
		\label{fig:features_q-dep}
	\end{figure*}
	
	\subsection{$D$-dependency}
	
	The statistics and evolution of homology classes of PCD are strongly dependent on the dimension of the Euclidean space as the embedding dimension. The zeroth (filled blue diamond) and first (empty orange square) PD of a regular ($q=0$) embedded fBm with $H=0.5$ by time-delay $\tau=1000$ to 2-, 3- and 4-cube is shown in the lower panels of Fig.~\ref{fig:H-dep_PD}. The inset plots are the projections of the reconstructed PCD into the standard 2-dimensional planes $P_{ij} \equiv \{ c_{i} \hat{e}_{i} + c_{j} \hat{e}_{j} ~ | ~ c_{i},c_{j} \in \mathbb{R} \}$, where $\{ \hat{e}_{i} \}_{i=1}^{D}$ is the standard basis of the space.  The  Betti curves for $d=0$ and $d=1$ of various $D$ are also shown in the lower panels of Fig.~\ref{fig:betti_curve_H-dep}. Since the typical distance between state vectors increases by increasing embedding dimension for any given value of Hurst exponent, the homology groups evolve (appeared, disappeared and maximized) in higher value of threshold (see also the lower panels of  Fig.~\ref{fig:H-dep_PD} as a consistent representation). %(Because of this we put the maximum computation scale as $\epsilon_{\rm max}=0.2 \sqrt{D}$). 
	The $H$-dependency of the measures increases by $D$ as shown in Fig.~\ref{fig:features_D-dep}. This means that the topological differences between the distribution of state vectors for various $H$ are more significant in higher dimensions which can be captured by the homology classes. However, the sensitivity of $E_0$ for the fBm signal with $H\lesssim 0.5$ is less than other defined measures irrespective of the embedding dimension.  
	
	The interesting thing about the $D$-dependency of the pattern of state vectors is that the topological distribution of these vectors for signals with $H\gtrsim0.5$ are more robust versus the embedding dimension $D$. It is worth noting that, the $D$-dependency of the homology groups evolution is considerable when the embedded time series is irregular (we will discuss in \ref{q-dependency}). 

	%^So we can write 
	%\begin{widetext}
	%\begin{equation}
	%	\frac{\partial \epsilon^{(0)}(H,D,\tau,q)}{\partial D}, \frac{\partial \epsilon_{\rm birth}^{(1)}(H,D,\tau,q)}{\partial D}, \frac{\partial \epsilon_{\rm death}^{(1)}(H,D,\tau,q)}{\partial D}, \frac{\partial \epsilon_{\rm max}^{(1)}(H,D,\tau,q)}{\partial D} \ge 0
	%\end{equation}
	%\end{widetext}
	%and for integration of the zeroth Betti curve 
	%\begin{equation}
	%	\frac{\partial B_{0}(H,D,\tau,q)}{\partial D} \ge 0
	%\end{equation}
	
	%Although, the width of the first Betti curve is an increasing function of $D$, but its height is effected by all parameters. 

	\subsection{$\tau$-dependency}
	
	%As a prior information from correlation analysis of time-series, the auto-correlation function of a fBm series decays in the power-law manner $c(\tau) \propto \tau^{-\alpha}$ $(\alpha=)$. This means that for such processes the behavior of the system become uncorrelated for long time lags. Therefore we expect that by increasing the time delay parameter $\tau$, $D$-dimensional state vectors become more scattered in unit $D$-cube. 
	We consider a  special case $\tau$=0 for which the embedded PCD, $X(x,D,\tau=0)=\{\vec x_{t} = (x_{t},...,x_{t})\}_{t=1}^{L}$, is a subset of a line along $\vec e \equiv \sum_{i=1}^{D} \hat{e}_{i}$ with trivial topology, $\tilde\beta_{d}(H,q,N,D,\tau=0,\eta)=0$, $d>0$, where $\{\hat{e}_{i}\}_{i=1}^{D}$ is the standard basis of $D$-dimensional Euclidean space. The upper panel of Fig.~\ref{fig:tau-dep_PD} illustrates 2-dimensional PCD reconstructed from highly correlated ($H=0.9$) fBm for various time-delay values $\tau=10, 100, 1000$, blue, red and orange, respectively (inset plots) and associated PDs (filled diamond for $d=0$ and empty square for $d=1$). The higher value of $\tau$ leads to construct more extended PCD as shown in the insets plots of upper panel of Fig.~\ref{fig:tau-dep_PD} and consequently, the amount of loops on the higher thresholds grows.

	The Betti curves of 10-dimensional PCD of fBm with $H=0.1$ for various value of $\tau$ also are shown in the upper panel of Fig.~\ref{fig:betti_curve_tau-dep}. Figure ~\ref{fig:features_tau-dep} reveals the behavior of our proposed topological measures as a function of the Hurst exponent for various values of $\tau$. Since the rate of the autocorrelation function decaying for the fBm time series decreases by increasing $H$ suggesting the small values of $\tau$ are more proper for estimating the Hurst exponent of the fBm series especially for $H\lesssim 0.5$ regime which is consistent with the statement in ~\cite{zou2015analyzing} from the network analysis point of view. For almost $H\lesssim 0.5$ regime, the autocorrelation of the fBm signal goes down rapidly compared to that for $H\gtrsim0.5$ and therefore as represented in different panels of Fig.~\ref{fig:features_tau-dep},  we find that for small $\tau$ our measures are as sensitive as enough to estimate Hurst exponent for the whole range $H\in(0,1)$ of underlying fBm signal. In addition, our results demonstrate that  particularly the $E_0$, $E_1$ and $n_1$ for large enough $\tau$ and small value of Hurst exponent loose their $H$-dependency.

	% The Fig.~\ref{fig:features_tau-dep} reveals that the behavior of some of our proposed measures looses its monotonic manner versus the Hurst exponent for large value of time-delay parameter in anti-correlated regime ($H<0.5$). The reason is that . This plot also indicates that the behavior of some features saturates to a distinct curve when we choose $\tau$ large enough. As we expect, by increasing $\tau$ the state vectors become more scattered in $\mathbb{R}^{D}$, so the critical scales shift to higher values. 
	
	%Therefore, increasing $\tau$ shifts the Betti curves to the high value of proximity parameter. But the main point is that for large enough value $\tau>>1$, the Betti curves do not move to the high thresholds anymore, Fig.\ref{fig:tau-dep_curve}

	\subsection{$q$-dependency}\label{q-dependency}

	The lower panel of Fig.~\ref{fig:tau-dep_PD} shows the zeroth (filled diamond) and first (empty square) PD of 2-dimensional PCDs (inset plots) converted from irregular ($q=0$ orange, $q=0.05$ red and $q=0.09$ blue) fBms with $H=0.1$ by $\tau=10$. The zeroth (middle left panel) and first (middle right panel) Betti curve of 2-dimensional PCDs mapped from the fBm signal with $H=0.1$ imposed by various irregularity $q$ are illustrated in Fig. (\ref{fig:betti_curve_tau-dep}). The lower panels of Fig. (\ref{fig:betti_curve_tau-dep}) are for 10-dimensional PCDs and the rest  parameters are the same as that for middle panels. This plot reveals the robustness of topological measures  when we consider low-dimensional PCDs against the irregularity in the time series. To explain this fact, by imposing irregularity to the time series $x$ of length $T$ makes the time series and  embedded PCD $X$ loose $T_{\rm missing}(q) = qT$ data points and $N_{\rm missing}(D,q) \ge T_{\rm missing}$ state vectors, respectively. The number of missing state vectors $N_{\rm missing}$ strongly depends on the embedding dimension, such that $N_{\rm missing}(D=1,q) = T_{\rm missing}$ and $N_{\rm missing}$ increases by $D$. Precisely, for a given irregular time series with irregularity equates to $q$, the number of missing data points in the series is $T_{\rm missing}=qT$, and according to $N=T-(D-1)\tau$, the influence of missing data point in time series grows for higher value of $D$ in state vectors of $X$.  Mentioned effect can be recognized in the zeroth Betti number of $X$  (see the beginning of zeroth Betti curve in the middle and lower left panels of  Fig.~\ref{fig:betti_curve_tau-dep}) as $N\tilde\beta_{0}(\eta=0) = N(q) \le T-(D-1)\tau$. This phenomenon affects the behavior of topological measurements versus the Hurst exponent. Such that some measures increase by $q$ but some others decrease. This effect dramatically is high in higher dimensions (see Fig.~\ref{fig:features_q-dep}).

	%On the other hand, typical distance between state vectors increases by $q$ which causes increasing path-connected regime critical scale for large $q$ , see Fig.~\ref{} 
	%\begin{equation}
	%	\frac{\partial \epsilon^{(0)}(H,D,\tau,q)}{\partial q} \ge 0
	%\end{equation}
	%Making the time series irregular also effects the first Betti curve. Statistically, increasing $q$ shifts the critical scales corresponding to the first homology groups to the higher values of proximity parameter 
	%\begin{equation}
	%	\frac{\partial \epsilon_{\rm birth}^{(1)}(H,D,\tau,q)}{\partial q}, \frac{\partial \epsilon_{\rm death}^{(1)}(H,D,\tau,q)}{\rm{d}q}, \frac{\partial \epsilon_{\rm max}^{(1)}(H,D,\tau,q)}{\partial q} \ge 0
	%\end{equation}
	%and decreases the height of the first Betti curve (Fig.~\ref{fig:q-dep_curve}(right)) 
	%\begin{equation}
	%	\frac{\partial \beta_{1}^{(max)}(\epsilon,H,D,\tau,q)}{\partial q}\le 0
	%\end{equation}
	\iffalse
	\begin{figure}
		\includegraphics[width=0.35\textwidth]{fig13.pdf}
		\caption{PDs of 2-dimensional PCDs (insets) for $q=0$ (orange), $q=0.05$ (red) and $q=0.09$ (blue).}
		\label{fig:q-dep_PD}
	\end{figure}
	\fi

	\section{Summary and Concluding remarks} \label{summary}

	In this work, we studied the topological signatures based on the homology groups of a point cloud data (PCD) constructed from the synthetic fractional Brownian motion (fBm) series to classify these series using the corresponding Hurst exponent.  We simulated the mock fBms for different $H$ and converted them to a $D$-dimensional PCD (a discrete subset of unit $D$-cube) according to the time delay embedding (TDE) method for various embedding dimensions and time delays. Then, by using the Vietoris-Rips (VR) method the embedded PCD is mapped into the simplicial complex for continuously varying proximity parameter $\epsilon$ and filtered by the persistent homology (PH) technique to capture the evolution of homology groups of triangulated PCD through the scales. For the zeroth and first homology group the homology classes are stored as pairs, so-called persistence pairs (PPs), in the zeroth and first persistence diagram (PD). Any PP contains appear and disappear scale of the homology class. 
	
	Relying on the population and distribution of PPs we defined the number of PPs and persistent entropy revealing the global properties of the corresponding PCD. Also, the zeroth and first Betti curves, $\beta_d(\epsilon)$, were computed to determine some transition scales for connectivity and loop structures in PCD.
	As we expected, all topological measures depend on the parameters ($H,q,D,\tau,\epsilon$) explained in our pipeline (subsection \ref{pipeline}). We assessed the $H$-dependency of topological measures considered in this paper. 
	
	Our results demonstrated that the $H$-dependency of our measures ($\eta_0^{\rm disappear}, B_0, E_0$ for zeroth homology group, $\eta_1^{\rm appear}, \eta_1^{\rm disappear}, \eta_1^{\rm maximize}, B_1, E_1, n_1$ for first homology group) grows by increasing  the embedding dimension (Fig.~\ref{fig:features_D-dep}).  The $D$-dependency goes down for the higher value of the Hurst exponent since for this range of $H$, the amount of autocorrelation becomes high for a small time-delay.   The time-delay imprint on the topological measures has been illustrated in Fig.~\ref{fig:features_tau-dep} demonstrating that the small value of $\tau$ is more reliable for determining the corresponding Hurst exponent in the whole range $H \in (0,1)$ of synthetic fBm.
	
	Motivated by irregularity existed in the realistic dataset in the universality class of fBm or fGn (\cite{ma2010effect,eghdami2018multifractal}), irregular mock series was simulated and quantified by a single irregular quantity ($q$). For the higher value of $H$, we expected the statistical properties of irregular fBm series remain almost unchanged. Our analysis showed that almost all topological criteria have weak dependency on $q$, meanwhile the $E_0$ ,  $E_1$ and $n_1$ depended on irregularity existed in data (Fig. \ref{fig:features_q-dep}). This means that those measures which are more related to the size and ordering of data are more sensitive to $q$, consequently, we propose that taking into account $\eta_0^{\rm disappear}$, $\eta_1^{\rm disappear}$  are more robust concerning data loss.  In addition, the $q$-dependency increases by increasing the embedding dimension.  Also, we find that selecting the lower value of time-delay makes a proper measure to estimate the value of the Hurst exponent for various types of fBm irrespective of irregularity.  
	
	For regular fBm series, the higher value of $D$, the more sensitive behavior of topological measures is. While the influence of irregularity would be magnified for the higher value of the embedding dimension yielding the more $H$-dependency happens for lower $D$.     
	
	The final remarks are as follows:  the size dependency of our proposed measures for any given value of Hurst exponent is important to examine.  Incorporating the embedding approach enables us to evaluate the higher dimensional topological holes which are nontrivial hidden shapes in the underlying series, particularly in the absence of irregularity. %\sadegh{However the latter needs more powerful computational hardware resources which are hardly available for us to complete nowadays}.  
	 Both mentioned tasks will be left for our next research.    
		
\section*{Acknowledgment}
The authors are very grateful to B. Askari for his constructive comments. A part of the numerical simulations were carried out on the computing cluster of the Brock University.

	%The situation is inverse when we deal with an irregular signal ($q>>0$). In this case, the proposed measures loose their robustness by increasing the embedding dimension (Figs.~\ref{fig:betti_curve_q-dep}~and~\ref{fig:features_q-dep}). 
	
	\iffalse
	\begin{figure}
		\includegraphics[width=0.49\textwidth]{fig141.pdf}
		\includegraphics[width=0.49\textwidth]{fig142.pdf}
		\caption{The number of connected components (left panels) and topological loops (right panels) of 2- (top panels) and 10-dimensional (bottom panels) PCD converted from fBm series ($H=0.1$) as a function of scale $\epsilon$ for various value of irregularity $q$. The behavior of topological features of the lower dimensional PCDs has more robust than higher dimensional ones in the presence of irregularity in the original signal.}
		\label{fig:betti_curve_q-dep}
	\end{figure}
	\fi

	\newpage
	
%	\bibliography{refs}

%	
\end{document}